\documentclass[12pt]{amsart}

\calclayout
\usepackage{amsmath,amsthm,amssymb}

\usepackage{colonequals}
\usepackage{enumitem,todo}

\usepackage[colorlinks=true,allcolors=blue,pdftitle={Unbounded mass solutions for the Keller-Segel equation}]{hyperref}

\usepackage[initials,nobysame,alphabetic,msc-links,backrefs]{amsrefs}%y2k

\numberwithin{equation}{section}

\newtheorem*{theorem*}{Theorem}
\newtheorem{thm}{Theorem}[section]

\newtheorem{prop}{Proposition}[section]
\newtheorem{lem}{Lemma}[section]

\newcommand{\la}{\lambda}
\newcommand{\ep}{\varepsilon}
\newcommand{\R}{\mathbb{R}}

\DeclareMathOperator{\diver}{div}
\newcommand{\rad}{\mathrm{rad}}

\title[The Keller--Segel equation in the disk]{Unbounded mass radial solutions for the Keller--Segel equation in the disk}

\author{Denis Bonheure}
\address{D\' epartement de Math\' ematique, Universit\' e libre de Bruxelles, Campus de la Plaine CP 213, Bd. du Triomphe, 1050 Bruxelles, Belgium} \email{denis.bonheure@ulb.ac.be}
\author{Jean-Baptiste Casteras}
\address{CMAFCIO, Faculdade de Ci\^encias da Universidade de Lisboa, Edificio C6, Piso 1, Campo Grande 1749-016 Lisboa, Portugal}
 \email{jeanbaptiste.casteras@gmail.com}
\author{Carlos Rom\'{a}n}
\date{June 1, 2021}
\address{Facultad de Matem\'aticas, Pontificia Universidad Cat\'olica de Chile, Vicu\~na Mackenna 4860, 7820436 Macul, Santiago, Chile}
\email{carlos.roman@mat.uc.cl}

\begin{document}
\begin{abstract}
We consider the boundary value problem
$$
\left\{ 
\begin{array}{rcll}
-\Delta u+ u -\la e^u&=&0,\ u>0 & \mathrm{in}\ B_1(0)\\
\partial_\nu u&=&0&\mathrm{on}\ \partial B_1(0),
\end{array}\right.
$$
whose solutions correspond to steady states of the Keller--Segel system for chemotaxis. Here $B_1(0)$ is the unit disk, $\nu$ the outer normal to $\partial B_1(0)$, and $\la>0$ is a parameter. We show that, provided $\la$ is sufficiently small, there exists a family of radial solutions $u_\la$ to this system which blow up at the origin and concentrate on $\partial B_1(0)$, as $\la\to 0$. These solutions satisfy 
$$
\lim_{\la\to 0} \frac{u_\la(0)}{|\ln\la|}=0\quad \mbox{and}\quad 0<\lim_{\la\to 0} \frac{1}{|\ln\la|}\int_{B_1(0)}\la e^{u_\la(x)}dx<\infty,
$$
having in particular unbounded mass, as $\la\to 0$.
\end{abstract}

\maketitle 

\section{Introduction}
Chemotaxis is the influence of chemical substances in the environment on the movement of mobile species. 
It is an important mean for cellular communication by chemical substances, which determines how cells arrange themselves, for instance in living tissues.
In 1970, Keller and Segel \cite{KelSeg} proposed a basic model to describe this phenomenon. They considered an advection-diffusion system consisting of two coupled
parabolic equations for the concentration of the species and that of the chemical released, respectively represented by strictly positive quantities
$v(x,t)$ and $u(x,t)$ defined on a bounded smooth domain $\Omega\subset \R^n$. This system has the form 
$$
\left\{ 
\begin{array}{rcll}
\dfrac{\partial v}{\partial t}&=&D_v\Delta v-c\diver (v\nabla \phi(u))&\mbox{in }\Omega\\[3mm]
\dfrac{\partial u}{\partial t}&=&D_u\Delta u+k(u,v)&\mbox{in }\Omega,\\[3mm]
\end{array}\right.
$$
with no flux through the boundary, that is, letting $\nu$ denote the exterior unit normal vector to $\partial\Omega$,
$$
\nabla v\cdot \nu =\nabla u\cdot \nu =0\quad \mathrm{on}\ \partial\Omega.
$$
Here, $D_v,D_u$, and $c$ are strictly positive constants, the function $\phi$, usually called the \emph{sensitive function}, is a smooth function such that $\phi'(u)>0$ for $u>0$ and $k$ is a smooth function such that $\frac{\partial k}{\partial v}\geq 0$ and $\frac{\partial k}{\partial u}\leq 0$. The typical choice for $k$ that we adopt from now on is $k(u,v)=-u+v$. 

An important property of this system is the so-called chemotactic collapse. This term refers to the fact that the whole population of organisms concentrate at a single point in finite or infinite time. When $\phi (u)=u$, it is well-known that the chemotactic collapse depends strongly on the dimension of the space. Finite-time blow-up never occurs if $n=1$, whereas it always occurs if $n\geq 3$. The two-dimensional case is critical: if the initial distribution of organisms exceeds a certain threshold, then the solutions may blow-up in finite time, whereas solutions exist globally in time if the initial mass is below the threshold. 
We refer the interested reader to the surveys \cites{Hor1,Hor2,MR3351175} and to the references therein for further details about the model and a collection of known results. 

Steady states of the Keller--Segel system are of basic importance for the understanding of the global dynamics. They solve the system
$$
\left\{ 
\begin{array}{rclcl}
-D_v \Delta v +c \diver (v \nabla \phi (u)) &=&0,&v>0&\mathrm{in}\ \Omega\\
-D_u \Delta u -u+v&=&0,&u>0&\mathrm{in}\ \Omega,
\end{array}\right.
$$
with homogeneous Neumann boundary conditions on $\partial\Omega$. This system can be reduced to a scalar equation as, 
indeed, one easily checks that
$$
\int_\Omega v |\nabla(D_v\ln v-c\phi(u))|^2\, dx=0,
$$
which implies $v=Ce^{\frac{c}{D_v}\phi(u)}$ for some constant $C>0$.
In the most common formulation of the Keller--Segel model, one takes $\phi(u)=u$, which yields the so-called Keller--Segel equation
\begin{equation}\label{ks}
\left\{ 
\begin{array}{rcll}
-\sigma^2 \Delta u+ u -\la e^u&=&0,\ u>0 & \mathrm{in}\ \Omega\\
\partial_\nu u&=&0&\mathrm{on}\ \partial\Omega,
\end{array}\right.
\end{equation}
where the constants $\sigma,\la$ depend on $D_v,D_u,c$ and $C$. 
It is worth mentioning that in the case $\phi(u)=\ln u$, one gets
$$
\left\{ 
\begin{array}{rcll}
-\tilde \sigma^2 \Delta u+ u -u^p&=&0,\ u>0 & \mathrm{in}\ \Omega\\
\partial_\nu u&=&0&\mathrm{on}\ \partial\Omega
\end{array}\right.
$$
for some constants $\tilde \sigma, p>0$, that is, one recovers the celebrated Lin--Ni--Takagi equation \cites{MR849484,MR974610,MR929196}. Let us observe that in dimension $2$ the Keller--Segel equation is critical, whereas the Lin--Ni--Takagi problem is subcritical. A good account of known results about this equation is given in the book by Wei and Winter \cite{WeiWinBook}, in the chapter \cite{MR2103689}, in the recent paper \cite{DelMusRomWei}, and in the references therein.

From now on, we study the Keller--Segel equation \eqref{ks} and we assume without loss of generality that $\sigma=1$. In the one-dimensional case, Schaaf \cite{schaaf} proved the existence of non-trivial solutions. For a general two-dimensional domain, the first existence results were obtained by  Wang and Wei \cite{WangWei2002} and independently by Senba and Suzuki \cite{MR1769174}, when the parameter $\lambda$ is small enough. Moreover, Senba and Suzuki \cites{MR1769174,MR1909263} studied the asymptotic behavior of finite mass solutions when $\lambda \rightarrow 0$. These are solutions $u_\lambda$ to \eqref{ks} such that
$$
\lim_{\lambda \rightarrow 0} \lambda \int_\Omega e^{u_\lambda}<\infty.
$$
They showed that there exist points $\xi_1,\dots,\xi_k \in \Omega$ and points $\xi_{k+1},\dots,\xi_{k+m} \in \partial \Omega$ such that, in the sense of measures,
\begin{equation}
\label{asymp}
-\Delta u_\la+u_\la=\la e^{u_\la}\rightharpoonup \sum_{i=1}^k 8\pi \delta_{\xi_i}+\sum_{i=k+1}^{k+m} 4\pi \delta_{\eta_i}
\end{equation}
in the sense of measures, and 
\begin{equation}
\label{profasy}
u_\lambda (x)\rightarrow \sum_{i=1}^k 8\pi \mathcal G(x,\xi_i)+\sum_{i=k+1}^m 4\pi \mathcal G(x,\eta_i),\quad \mathrm{as}\ \lambda \to 0,
\end{equation}
uniformly on compact subsets of $\overline{\Omega}\backslash \{\xi_1,\ldots ,\xi_k,\eta_{k+1},\dots,\eta_{m}\}$, where, given $y\in \overline \Omega$, $\mathcal G(x,y)$ denotes the Green function that uniquely solves
$$
\left\{ 
\begin{array}{rcll}
-\Delta_x \mathcal G+\mathcal G&=&\delta_y & \mathrm{in}\ \Omega\\
\nabla \mathcal G\cdot \nu&=&0&\mathrm{on}\ \partial\Omega.
\end{array}\right.
$$

The counterpart of this result was obtained by del Pino and Wei \cite{MR2209293}. For any given integers $k$ and $m$, they constructed a family of solutions $u_\la$ to \eqref{ks} that satisfy \eqref{asymp} and \eqref{profasy} for a suitable choice of points $\xi_i\in \Omega$ for $i=1,\dots,k$ and $\xi_i\in \partial\Omega$ for $i=k+1,\dots,k+m$. Near each of these points $\xi=\xi_i$,
$$
u_\la(x)\approx V_{\mu_i}(|x-\xi|),
$$
where $V_{\mu_i}$ is a radially symmetric solution to
$$
-\Delta V-\la e^V=0\quad \mbox{in }\R^2,
$$
that is, a function of the form
\begin{equation}\label{near0}
V_{\mu}(|x|)=\ln \frac{8\mu^2}{(\la\mu^2+|x|^2)^2}\quad \mbox{with }\mu>0.
\end{equation}
The parameter $\mu_i=C(\xi_1,\dots,\xi_{k+m},\Omega)$ is a constant that depends only on the points $\xi_i$'s and $\Omega$. In particular, it does not depend on $\lambda$. 

\medskip
It is worth mentioning that these solutions have quantized mass, that is,
$$
\lim_{\la\to 0}\int_\Omega \la e^{u_\la}=4\pi(2k+m).
$$

\medskip
Recently, solutions concentrating on higher dimensional sets with unbounded mass, namely
$$
\lim_{\la\to 0}\int_\Omega \la e^{u_\la}=\infty,
$$
have been proven to exist. From now on, we denote by $B_r$ the ball of radius $r$ centered at zero.
When $\Omega=B_1\subset \R^n$ with $n\geq 2$, Pistoia and Vaira \cite{PistoiaVaira2015} constructed a family $u_\lambda$ of radial solutions concentrating on the whole boundary of $\Omega$ such that
$$0<\lim_{\lambda \rightarrow 0} \frac{1}{|\ln \la|} \int_{B_1} \la e^{u_\lambda (x)}dx<\infty $$
More precisely, their solutions satisfy 
$$\lim_{\lambda \rightarrow 0}\frac1{|\ln \la|} u_\lambda=\mathcal U,$$
$C^0$-uniformly on compact subsets of $B_1$, where 
$\mathcal U$ is the unique (radial) solution to
\begin{equation*}
\left\{ 
\begin{array}{rcll}
-\Delta \mathcal U +  
\mathcal U&=&0& \mathrm{in}\ B_1\\
\mathcal U&=&1&\mathrm{on}\ \partial B_1.
\end{array}\right.
\end{equation*}
whereas near the boundary,
\begin{equation}\label{near1}
u_\la(r) +\ln \lambda \approx W_{\varepsilon}(r)= \ln \left( \frac 4{\varepsilon^2}\frac{e^{\sqrt2\frac{r-1}{\varepsilon}}}{\left(1+e^{\sqrt2\frac{r-1}{\varepsilon}}\right)^2}\right),
\end{equation}
where $\ep=\ep_\la\approx \frac{\sqrt{2}}{\mathcal U'(1)}\frac{1}{|\ln \lambda|}$.
 
Let us point out that $W_{\ep}$ is a radial solution to the one-dimensional problem 
$$
-W^{\prime \prime}=e^W\quad \text{in}\ \R, \quad \text{with}\ \int_{\R} e^W<\infty.
$$

\smallskip
Del Pino, Pistoia, and Vaira \cite{del2014large}, generalized this result to general two-dimensional domains. Also, the existence of solutions concentrating on submanifolds of the boundary has also been investigated; see for instance \cite{agpis}.
 
\medskip
From now on, we suppose that $\Omega=B_1\subset \R^n$, with $n\geq 2$. In \cite{bocano}, a bifurcation analysis of radial solutions to \eqref{ks} was performed. Observe that for $\lambda <1/e$, the equation \eqref{ks} can be rewritten as
\begin{equation}\label{eq:intro_u_mu}
\left\{
\begin{array}{rcll}
-\Delta u+u&=&e^{\mu(u-1)},\ u>0 &\text{in } B_1\\
\nabla u\cdot \nu &=&0&\text{on } \partial B_1
\end{array}\right.
\end{equation}
for $\mu>1$. This equation admits the constant solutions $u\equiv 1$ and  $\underline{u}_\mu <1$. To describe the bifurcation result, we denote by $\lambda_i^{\text{rad}}$ the $i$-th eigenvalue of the operator $-\Delta+\text{Id}$ in $B_1$, restricted to the set of radial functions, with homogeneous Neumann boundary conditions.

\begin{theorem*}[\cite{bocano}]
For every $i\geq 2$, $(\lambda_i^{\rad},1)$ is a bifurcation point of \eqref{eq:intro_u_mu}. Denoting by $\mathcal{B}_i$ the continuum that branches out of $(\lambda_i^{\rad},1)$, we have
\begin{enumerate}[leftmargin=*,label=(\roman*)]
\item the branches $\mathcal{B}_i$ are unbounded and do not intersect; close to $(\lambda_i^{\rad},1)$, $\mathcal{B}_i$ is a $C^1$-curve;
\item if $u_\mu \in \mathcal{B}_i$ then $u_\mu >0$;
\item each branch consists of two connected components: the component $\mathcal{B}_i^-$, along which $u_\mu(0)<1$, and the component $\mathcal{B}_i^+$, along which $u_\mu(0)>1$;
\item if $u_\mu \in \mathcal{B}_i$ then $u_\mu-1$ has exactly $i-1$ zeros, $u_\mu'$ has exactly $i-2$ zeros, and each zero of $u_\mu'$ lies between two zeros of $u_\mu-1$;
\item the functions satisfying $u_\mu(0)<1$ are uniformly bounded in the $C^1$-norm.
\end{enumerate}
\end{theorem*}

We conjecture that the solutions constructed by Pistoia and Vaira \cite{PistoiaVaira2015} correspond to those on $\mathcal{B}_1^-$, while the solutions constructed by del Pino and Wei \cite{MR2209293} (when restricted to the $2$-dimensional ball) correspond to the branch $\mathcal{B}_1^+$. In \cite{bocano}, the authors constructed solutions concentrating on an arbitrary number of internal spheres by combining variational and perturbative methods. Solutions sharing the same qualitative properties were obtained with a different method in \cite{bocano2} with very precise asymptotics. We conjecture that those solutions are indeed the same and correspond to the solutions on the branches $\mathcal{B}_i^-$. 

In this paper, we restrict ourselves to the disk, that is, the case $n=2$. Our main result is the construction of solutions to \eqref{ks} that concentrate at the origin and on the boundary of $B_1$, as $\lambda\to 0$.
We conjecture that they correspond to the solutions to \eqref{eq:intro_u_mu} on the branch $\mathcal{B}_2^+$. We emphasize that only a few results concerning existence of solutions concentrating simultaneously on points and layers are available in the literature, see for instance \cites{MR3038717,MR2379470}. 

\begin{thm}\label{mainthm} 
There exist $\la_0>0$ and a family of radial solutions $\{u_\lambda\mid\,  \lambda \in (0,\lambda_0)\}$
to \eqref{ks} such that
$$
\lim_{\la\to 0} \frac{u_\la(0)}{|\ln \lambda |}=\infty \quad \mbox{and}\quad 0<\lim_{\lambda \rightarrow 0} \frac1{|\ln\lambda|}\int_{B_{1}}\lambda e^{u_\lambda (x)} dx <\infty.
$$
Moreover, letting $\varepsilon_\lambda\to 0$ as $\lambda\to 0$ be the parameter defined by 
\begin{equation}\label{epla}
\ln \dfrac{4}{\ep_\la^2}-\ln \lambda = \dfrac{a_{1,\ep_\la}}{\ep_\la}+a_{2,\ep_\la} +a_{3,\ep_\la} \ep_\la
\end{equation}
for suitable constants $a_{1,\ep_\la}$, $a_{2,\ep_\la}$, and $a_{3,\ep_\la}$ depending on $\ep_\la$ (see \eqref{ai}), and 
letting $G_{\ep_\la}$ be the unique radial solution for $\ep=\ep_\la$ to
\begin{equation}
\label{greenintro}
\left\{
\begin{array}{rcll} 
-\Delta G_{\ep}
+G_{\ep}&=&0 &\mathrm{in }\ B_1\\[2mm] 
\displaystyle\lim_{r\rightarrow 0^+} \dfrac{G_{\ep} (r)}{|\ln r|}&=&4/{\mathcal A_{\ep}}&\\
G_{\ep} &=&1 & \mathrm{on }\ \partial B_1,\end{array}\right.
\end{equation}
where $\mathcal{A}_{\ep}=\frac{A_{1,\ep}}{\ep} + A_{2,\ep} +A_{3,\ep}\ep$ for suitable constants $A_{1,\ep},A_{2,\ep},$ and $A_{3,\ep}$ depending on $\ep$ (see \eqref{Ai}), we have that, uniformly on compact subsets of  $B_{1}\backslash \{0\}$,
$$
\lim_{\lambda \rightarrow 0} \left( u_\lambda(x) -  \mathcal{A}_{\ep_\la} G_{\ep_\la}(x) \right)=0,
$$ 
and, in the sense of measures,
$$
\lambda e^{u_\lambda} \rightharpoonup 8\pi \delta_{0}\quad \mathrm{in}\ B_{1/2},\quad\quad \ep_\la \la e^{u_\lambda} \rightharpoonup \sqrt2 \delta_{\partial B_1} \quad \mathrm{in}\ \overline B_{1}\backslash  B_{1/2}.
$$
\end{thm}

Theorem \ref{mainthm} will be proven using a fixed point argument. More precisely, we will look for a solution to \eqref{ks} of the form $u_\la=U_\la+\phi_\la$, where $U_\la$ is a first ``good'' approximation of the solution and $\phi_\la$ is a small perturbation. Roughly speaking, $U_\la$ is constructed by gluing the Green's function $\mathcal A_{\ep_\la} G_{\ep_\la}$ with $V_{\mu}$ (recall \eqref{near0}) near the origin and with $W_{\ep}$ (recall \eqref{near1}) near the boundary, for well chosen parameters $\mu$ and $\ep$. To obtain a ``good'' matching between these functions, we are forced to choose $\varepsilon_\la$ satisfying \eqref{epla} and 
$$
\mu^2=\mu_\la^2 = \frac{e^{H_{\ep_\la}(0)}}8,
$$ 
where $H_{\ep_\la}$ denotes the regular part of $\mathcal A_{\ep_\la} G_{\ep_\la}$, that is $H_{\ep_\la}(r)=\mathcal A_{\ep_\la} G_{\ep_\la}(r)+4\ln r$. We explicitly compute $H_{\ep_\la}(0)$ (see \eqref{estrobin}), which leads to 
$$
8\mu_\la^2\approx e^\frac{A_1}{C\ep_\la}\to \infty \quad\mbox{as }\la\to 0, 
$$
where $A_1=\lim_{\ep \to 0}A_{1,\ep}>0$ and $C>1$ is a constant. Let us stress that this is very different to the previously described situation of finite mass blow-up, and poses technical issues. In fact, what allows our argument to work is the crucial fact that
$$
\lim_{\la\to 0}\la \mu_\la^2=\lim_{\la \to 0} \frac{1}{\ep_\la^2}\exp\left(\frac{A_1}{\ep_\la}\left(\frac1C-1\right)\right)=0.
$$

\medskip
The rest of the paper is organized as follows. In Section \ref{sec2}, we provide the existence of the Green's function solution to \eqref{greenintro}, which is used to build the first approximation of the solution in Section \ref{sec3}. We then estimate the error introduced by our approximation in Section \ref{sec4}. In Section \ref{sec5}, we prove the solvability of the linearized equation around our approximate solution, which allows us to use a fixed point argument to prove Theorem \ref{mainthm} in Section \ref{secMain}.

\section*{Acknowledgments}
D. Bonheure is partially supported by the project ERC Advanced Grant  2013 n. 339958: ``Complex Patterns for Strongly Interacting Dynamical Systems - COMPAT'' and by ARC AUWB-2012-12/17-ULB1- IAPAS.

J-B. Casteras is supported by MIS F.4508.14 (FNRS), PDR T.1110.14F (FNRS) and INRIA Team MEPHYSTO.

This work was initiated when C. Rom\'an was a Ph.D. student
at the Jacques-Louis Lions Laboratory of the Pierre and Marie Curie University, supported by a public grant overseen by the French National Research Agency (ANR) as part of the ``Investissements d’Avenir" program (reference: ANR-10-LABX-0098, LabEx SMP). He is currently supported by the Chilean National Agency for Research and Development (ANID) through FONDECYT Iniciaci\'on grant 11190130. He wishes to thank the support and kind hospitality of the Universit\'e libre de Bruxelles, where part of this work was done.

\section{Green's function}\label{sec2}
This section is devoted to prove the existence of $G_\ep$ defined in \eqref{greenintro}. First, let us recall the following lemma from \cite{bocano}.
\begin{lem}
\label{lemma:appendix1}
There exist two positive linearly independent solutions $\zeta \in C^2 ((0,1])$ and $\xi \in C^2 ([0,1])$ of the modified zero-order Bessel differential equation
$$-u^{\prime \prime}-\dfrac{1}{r}u^\prime +u=0\quad \mathrm{in}\ (0,1),$$
satisfying
$$\xi^\prime (0)=\zeta^\prime (1)=0\quad \mathrm{and}\quad r(\xi^\prime (r) \zeta (r) - \xi (r) \zeta^\prime (r))=1\ \mathrm{for\ any}\ r\in (0,1].$$
We have that $\xi$ is bounded and strictly increasing in $[0,1]$, $\zeta$ is strictly decreasing in $(0,1]$, 
$$\xi (0)=1,\quad \lim_{r\rightarrow 0^+}\dfrac{\zeta (r)}{|\ln r|}=1,\quad \mathrm{and} \quad \lim_{r\rightarrow 0^+ }(-r\zeta^\prime (r))=1.$$
Moreover, as $r$ goes to $0$, we have (see \cite{bookBessel})
\begin{equation}
\label{besselasympzero}
\zeta (r)=  (|\ln r| +c_1) +\dfrac{r^2}{4} (|\ln r| +c_2)+O(r^4|\ln r|),
\end{equation}
\begin{equation}
\label{derbesselasympzero}
\zeta ' (r)= -\dfrac{1}{ r} + \dfrac{r(|\ln r| +c_3)}{2}+O(r^3|\ln r|),
\end{equation}
and 
\begin{equation}\label{bessel2}
\xi(r)=1+\frac{r^2}4+O(r^4),\quad \xi'(r)=\frac r2 +O(r^3),
\end{equation}
where $c_1,c_2,c_3$ are positive constants.
\end{lem}

Using this result, we are able to construct a radial Green's function on the unit ball $B_1$ blowing up at $0$ and equal to $1$ on $\partial B_1$.
\begin{lem}
\label{singreen}
For any $\varepsilon>0$ small enough, there exists a function $G_\varepsilon$ satisfying
\begin{equation}
\label{greensing}
\left\{
\begin{array}{rcll}
- G_\varepsilon^{\prime \prime}-\dfrac{1}{r}G_\varepsilon^\prime +G_\varepsilon&=&0&\mathrm{in}\ (0,1)\\
\displaystyle \lim_{r\rightarrow 0^+}\dfrac{G_\varepsilon (r)}{|\ln r|}&=&4/\mathcal{A}_\varepsilon &\\
G_\varepsilon (1)&=&1,&
\end{array}\right.
\end{equation}
where $\mathcal{A}_\varepsilon=\dfrac{A_{1,\varepsilon}}{\varepsilon} +A_{2,\varepsilon} +A_{3,\varepsilon} \varepsilon$, with
$$A_{1,\varepsilon}=\dfrac{\sqrt{2}}{G_\varepsilon^\prime (1)},\ A_{2,\varepsilon}=\dfrac{1}{G_\varepsilon^\prime (1)} \left(\dfrac{\ln 4}{G_\varepsilon^\prime (1)}-2\right),\ A_{3,\varepsilon}=\dfrac{c}{G_\varepsilon^\prime (1)}$$
for $c\in \R$.

Moreover, there exists $\tilde{r}\in (0,1)$ with $\tilde r\approx \sqrt{\varepsilon}$, that is, there exist two constants $c_1,c_2>0$ such that $c_1 \sqrt{\varepsilon}\leq \tilde r \leq c_2 \sqrt{\varepsilon}$, such that $G_\varepsilon^\prime (\tilde{r})=0$, and there holds $\lim_{\varepsilon \rightarrow 0}G_\varepsilon^\prime (1)=\dfrac{\xi^\prime (1)}{\xi (1)}$.
We also have, as $r$ goes to zero,
\begin{equation}
\label{besselasympzerobis}
G_\varepsilon (r)-\dfrac{4}{\mathcal{A}_\varepsilon}|\ln r|=\dfrac{1}{\xi (1)}+O_r(r^2)+\frac{O_r(r^2|\ln r|)}{\mathcal{A}_\ep}+O_{\ep}(\ep),
\end{equation}
where $C$ is a constant independent of $\ep$, and
\begin{equation}
\label{derbesselasympzerobis}
G_\varepsilon^\prime (r) + \dfrac{4}{\mathcal{A}_\varepsilon r}=O(r)+\frac{O(r|\ln r|)}{\mathcal A_\ep} .
\end{equation}
\end{lem}
\begin{proof}
Using the properties of the functions $\xi$ and $\zeta$ (defined in Lemma \ref{lemma:appendix1}), it is immediate to see that, for any $b\in (0,1)$,
$$u_b (r)=\dfrac{\xi^\prime (b) \zeta (r) -\xi (r) \zeta^\prime (b) }{\xi^\prime (b) \zeta (1) -\xi (1) \zeta^\prime (b)}$$
is a solution to \eqref{greensing} such that 
$$u_b (1)=1\quad \mbox{and}\quad \lim_{r\rightarrow 0^+} \dfrac{u_b (r)}{|\ln r|}=\dfrac{\xi^\prime (b)}{\xi^\prime (b) \zeta (1) -\xi (1) \zeta^\prime (b)}.$$
Moreover, for $b$ small enough, we have
$$
\begin{array}{rcl}
\xi^\prime (b) &=& \frac12 b +o(b),\\
\xi^\prime (b) \zeta (1)  - \xi (1) \zeta^\prime (b)&=& \xi (1) b^{-1}+O(b).
\end{array}
$$
Therefore, for $b$ small enough, we have
$$ \lim_{r\rightarrow 0^+} \dfrac{u_b (r)}{|\ln r|}= \frac1{2\xi(1)} b^2+o(b^2).$$
On the other hand, it is easy to check that 
\begin{equation}\label{li}
u_b^\prime (1)= \dfrac{\xi^\prime (1)}{\xi (1)}+o_b (1),
\end{equation}
where $o_b (1)\rightarrow 0$ as $b\rightarrow 0$. Therefore, we can choose $b\approx \sqrt\ep$ such that  $\displaystyle\lim_{r\rightarrow 0^+} \dfrac{u_b (r)}{|\ln r|}=4/\mathcal{A}_\varepsilon$, which proves the existence of the function $G_\ep$. 

From \eqref{li}, we immediately see that $\lim_{\ep\to 0}G_\ep '(1)=\frac{\xi'(1)}{\xi(1)}$. The fact that $u'_b(b)=0$ implies the existence of $\tilde r\approx \sqrt \ep$ such that $G_\ep'(\tilde r)=0$. Finally, \eqref{besselasympzerobis} and \eqref{derbesselasympzerobis} follow from \eqref{besselasympzero}, \eqref{derbesselasympzero}, and \eqref{bessel2}.
\end{proof}

\section{The approximate solution}\label{sec3}
We look for a radial solution to \eqref{ks} concentrating at $0$ and on $\partial B_1 $. To do so, we take an ansatz of solution of the form
$$
U=\left\{
\begin{array}{cl}u_0& \mbox{in}\ [0,\delta)\\u_1& \mbox{in}\ [\delta,2\delta)\\u_2& \mbox{in}\ [2\delta,1-2\delta_1)\\u_3&\mbox{in}\ [1-2\delta_1,1-\delta_1)\\u_4&\mbox{in}\ [1-\delta_1,1],
\end{array}\right.
$$
where $\delta$ and $\delta_1$ are suitable constants depending on $\lambda$. Let us first describe our ansatz in words. Near the origin, we want  $U=u_0$ to behave approximately like $V_\mu$, the two dimensional standard bubble given by
\begin{equation}
\label{defU_0}
V_\mu(r)= \ln \dfrac{8\mu^2}{(\la\mu^2+r^2 )^2},
\end{equation}
for some constant $\mu=\mu_\la>0$ to be specified later. Let us recall that these functions correspond to all solutions of the problem
$$
-\Delta V=\lambda e^V\quad \mbox{in }\R^2,\quad \mbox{with } \lambda \int_{\R^2}e^Vdx<\infty.
$$
Near the boundary $\partial B_1$ of the disk, we want that $U=u_4$ behaves like $W_{\ep} -\ln \la$ where $W_{\ep}$ is the one dimensional standard bubble solving $-W''=e^W$ in $\R$, which is given by
\begin{equation*}
W_{\varepsilon} (r) = \ln \left( \dfrac{4}{\varepsilon^2}\dfrac{e^{-\frac{\sqrt{2} (r-1)}{\varepsilon}}}{\left(1+e^{-\frac{\sqrt{2} (r-1)}{\varepsilon}}\right)^2} \right),
\end{equation*}
for some constant $\ep=\ep_\la>0$ to be determined later. In order to ``glue'' these singular solutions, far from the origin and $\partial B_1 $ we choose $U=\mathcal A_{\ep_\la} G_{\ep_\la}$, where $G_{\ep_\la}$ is the singular at the origin Green's function given in Lemma \ref{singreen} (with $\ep=\ep_\la$ and a suitable constant $c$). Finally, we choose $u_1$ and $u_3$ to be linear interpolations between $u_{i-1}$ and $u_{i+1}$, for $i=1,3$, namely 
\begin{equation}
\label{defu_1}
u_i (r)=\chi_i (r) u_{i-1}(r)+(1-\chi_i (r)) u_{i+1}(r),
\end{equation}
where $\chi_i \in C^2 ((0,1))$ are cut-off functions such that
$$\chi_1 (r)\equiv 1\ \mbox{in}\ (0,\delta),\ \chi_1 (r)\equiv 0\ \mbox{in}\ (2\delta ,1 ),\ |\chi_1 (r)|\leq 1,\ |\chi_1^\prime (r)|\leq c\delta^{-1} , |\chi_1^{\prime \prime}(r)|\leq c\delta^{-2},$$
and
$$\chi_3 \equiv 1\ \mbox{in}\ (0,1- 2\delta_1),\ \chi_3 \equiv 0\ \mbox{in}\ (1- \delta_1 ,1 ),\ |\chi_3 (r)|\leq 1,\ |\chi_3^\prime (r)|\leq c\delta_1^{-1} , |\chi_3^{\prime \prime}(r)|\leq c\delta_1^{-2}.$$

We now describe our ansatz in detail.

\subsection{Construction of \texorpdfstring{$u_4$}{u4}}
First, let us set $\ep=\ep_\la\to 0$ as $\la\to 0$, via the relation
\begin{equation}
\label{relpara}
\ln \dfrac{4}{\varepsilon^2}-\ln \lambda = \dfrac{a_{1,\ep}}{\varepsilon}+a_{2,\ep} +a_{3,\ep} \varepsilon,
\end{equation}
for some constants $a_{i,\ep}$, $i=1,2,3$ to be determined later (see \eqref{ai}), and let
\begin{equation}
\label{defdelta1}
\delta_1=\varepsilon^{\eta},\text{ for some }\eta \in \left(\frac23,1\right).
\end{equation}
We define $u_4$ in the same way as the function ``$u_1$'' of \cite{PistoiaVaira2015} (or \cite{bocano2}) with $r_0=1$. The construction of this function is quite lengthy .We only briefly recall it here, and refer to the above two papers for more details. We define
\begin{equation*}
u_{4}=\underbrace{W_{\varepsilon}-\ln \lambda+\alpha_\varepsilon}_{1^{st}\ order\ approx.} +\underbrace{v_\varepsilon+\beta_\varepsilon}_{2^{nd}\ order}+\underbrace{z_\varepsilon}_{3^{rd}\ order},
\end{equation*}
where $\alpha_\ep(r),v_\ep(r),\beta_\ep(r),$ and $z_\ep(r)$, which we briefly describe below, are functions introduced in order to produce a good enough match between $u_4$ and $u_2= \mathcal A_\ep G_\ep$; see Lemma \ref{lemrelu2u4}.
 
\medskip
The function $\alpha_\varepsilon$ satisfies
$$
\left\{
\begin{array}{rcll}
-(\alpha_\varepsilon )''-\dfrac{1}{r}(\alpha_\varepsilon )^\prime&=&\dfrac{1}{r}(W_\varepsilon )^\prime - W_\varepsilon +\ln \lambda&\mathrm{in}\ (0,1)\\
\alpha_\varepsilon(1)&=&0&\\
(\alpha_\varepsilon)'(1)&=&0&
\end{array}\right.
$$
and the following estimate holds, for $s\leq 0$,
\begin{equation}
\alpha_\varepsilon (\varepsilon s+1)=\varepsilon (\alpha_\varepsilon)_1 (s)+\varepsilon^2 (\alpha_\varepsilon)_2 (s)+O (\varepsilon^3 s^4),
\label{expbound}
\end{equation}
where, letting $W$ be defined via
\begin{equation*}
W\left(\dfrac{r-1}{\varepsilon}\right)+\ln \dfrac{4}{\varepsilon^2}-\ln 4=W_{\varepsilon}(r),
\end{equation*}
we have 
$$(\alpha_\varepsilon)_1 (s)= - \int_0^s W (\sigma )d\sigma +\dfrac{a_{1,\la}}{2} s^2$$
and
\begin{align*}
(\alpha_\varepsilon) _2 (s) & = \int_0^s \int_0^\sigma (W (\rho) -\ln 4)d\rho d\sigma + \int_0^s \sigma W (\sigma )d\sigma -\dfrac{1}{6} a_{1,\la} s^3 +\dfrac{a_{2,\la}}{2} s^2 .
\end{align*}

\medskip
The function $v_\varepsilon$ solves
$$
\left\{
\begin{array}{rcll}
-(v_\varepsilon)'' - e^{W_{\varepsilon}} v_\varepsilon&=&\varepsilon e^{W_{\varepsilon}} (\alpha_\varepsilon)_1 \left(\dfrac{r-1}{\varepsilon}\right)&\text{in}\ \R\\
v_\varepsilon(1) &=&0&\\
(v_\varepsilon)'(1)&=&0&,
\end{array}\right.
$$
where $(\alpha_\varepsilon)_1$ is defined in \eqref{expbound}. Moreover, we have
\begin{equation*}
v_\varepsilon (r)=\nu_1 (r-1)+\nu_2 \varepsilon+O(\varepsilon e^{-\frac{|r-1|}{\varepsilon}}),
\end{equation*}
where
\begin{equation*}
\nu_2 \in \R \quad \mathrm{and}\quad  \nu_1= -2(1-\ln 2) +a_1 \sqrt{2} \ln 2 .
\end{equation*}

\medskip
The function $\beta_\varepsilon$ is the solution of
$$
\left\{
\begin{array}{rcll}
-(\beta_\varepsilon)'' - \dfrac{1}{r}(\beta_\varepsilon)'&=&\dfrac{1}{r} (v_\varepsilon)'&\mathrm{in}\ (0,1),\\
\beta_\varepsilon (1)&=&0,&\\
(\beta_\varepsilon)'(1)&=&0,&
\end{array}\right.
$$
and the following estimate holds, for $s\leq 0$,
\begin{equation*}
\beta_\varepsilon (\varepsilon s+1)=\varepsilon^2 (\beta_\varepsilon)_1 (s) +O(\varepsilon^3 s^3),
\end{equation*}
where, letting $v$ be defined via
\begin{equation*}
v_\varepsilon (r)= \varepsilon v\left(\dfrac{r-1}{\varepsilon}\right),
\end{equation*}
we have
$$ (\beta_\varepsilon)_1 (s)= - \int_0^s \int_0^\sigma v'(\rho) d\rho d\sigma.$$

\medskip
Finally, the function $z_\varepsilon$ satisfies
\begin{equation*}
\left\{
\begin{array}{rcll}
-(z_\varepsilon)'' -e^{W_{\varepsilon}} z_\varepsilon&=&\varepsilon^2 e^{W_{\varepsilon}} \left[(\alpha_\varepsilon)_2 \left(\dfrac{r-1}{\varepsilon}\right) +(\beta_\varepsilon)_1  \left(\dfrac{r-1}{\varepsilon}\right)\right.&\\
&&
\quad\left.+\dfrac{1}{2} \left((\alpha_\varepsilon)_1  \left(\dfrac{r-1}{\varepsilon}\right) +v  \left(\dfrac{r-1}{\varepsilon}\right)\right)^2 \right]&\mathrm{in} \ (0,1)\\
z_\varepsilon (1)&=&0&\\
(z_\varepsilon )^\prime (1)&=&0
\end{array}\right.
\end{equation*}
and there holds
\begin{equation*}
z_\varepsilon(r) =\varepsilon \zeta_1 (r-1) +\zeta_2 \varepsilon^2 +O(\varepsilon^2 e^{-\frac{|r-1|}{\varepsilon}})
\end{equation*}
for some $\zeta_j \in \R$, $j=1,2$.

\subsection{Construction of \texorpdfstring{$u_2$}{u2}}
Thanks to Lemma \ref{singreen} (with $c=\zeta_1$), we know that, for any $\ep$ small enough, there exists a function $G_{\ep}$ satisfying
$$
\left\{
\begin{array}{rcll}
- G_\ep^{\prime \prime}-\dfrac{1}{r}G_\ep^\prime +G_\ep&=&0&\mathrm{in}\ (0,1)\\
\displaystyle \lim_{r\rightarrow 0^+}\dfrac{G_\ep(r)}{|\ln r|}&=&4/\mathcal{A}_\varepsilon &\\
G_\ep(1)&=&1.&
\end{array}\right.
$$
where $\mathcal{A}_\ep=\dfrac{A_{1,\ep}}{\ep} +A_{2,\ep} +A_{3,\ep} \ep$ and
\begin{equation}\label{Ai}
A_{1,\ep}=\dfrac{\sqrt{2}}{G_\ep^\prime (1)},\quad A_{2,\ep}=\dfrac{1}{G_\ep^\prime (1)} \left(\dfrac{\ln 4}{G_\ep^\prime (1)}-2\right),\quad A_{3,\ep}=\dfrac{\zeta_1}{G_\ep^\prime (1)}.
\end{equation}

Letting
\begin{equation}\label{ai}
a_{1,\ep}=A_{1,\ep},\quad a_{2,\ep}=A_{2,\ep},\quad a_{3,\ep}=A_{3,\ep}-\nu_2,
\end{equation}
and recalling \eqref{relpara}, we define
\begin{equation}
\label{defu_2}
u_2(r)=\mathcal{A}_{\ep_\la}  G_{\ep_\la} (r). 
\end{equation}
Thanks to our definition of $u_2$ and $u_4$, one can show, arguing as in \cite{bocano2}*{Lemma 3.3}, the following estimates.
\begin{lem}
\label{lemrelu2u4}
For any $\delta_1 <|r-1|<2\delta_1 $, we have
$$u_4 (r)-u_2 (r)= O\left(\ep_\la^2+\ep_\la |r-1|^2+|r-1|^3 + \dfrac{|r-1|^4}{\ep_\la}+\exp\left(-\dfrac{|r-1|}{\ep_\la}\right)\right)$$
and
$$u_4^\prime (r)-u_2^\prime (r)= O\left(\ep_\la |r-1|+|r-1|^2 + \dfrac{|r-1|^3}{\ep_\la}+\dfrac{1}{\ep_\la}\exp\left(-\dfrac{|r-1|}{\ep_\la}\right)\right).$$
\end{lem}

\medskip
In order to define $u_0$ and estimate $u_0-u_2$, it is important to introduce the regular part $H_{\ep_\la}$ of $u_2$, namely
\begin{equation}
\label{regreen}
H_{\ep_\la}(r)=u_2 (r) +4\ln r. 
\end{equation}
We let $\mu_\la>0$ be defined via the relation
\begin{equation}\label{defmu}
8\mu_\la^2=e^{H_{\ep_\la}(0)}.
\end{equation}
Thanks to \eqref{besselasympzerobis}, \eqref{derbesselasympzerobis}, and $\lim_{\la \rightarrow 0} G_{\ep_\la}^\prime (1)=\dfrac{\xi^\prime (1)}{\xi (1)}$, we have, for a constant $c$ independent of $\ep_\la$,
\begin{equation}\label{estrobin}
H_{\ep_\la}(0)=\dfrac{\sqrt{2}}{  \xi^\prime (1)\ep_\la}+c+O_{\ep_\la}(\ep_\la)\quad \text{and}\quad \lim_{r\rightarrow 0^+} H_{\ep_\la}'(r)=0.
\end{equation}

Moreover, as $r\to 0$, we have
\begin{equation}\label{expansionH}
|H_{\ep_\la}(r)-H_{\ep_\la}(0)|\leq C\left( \frac{r^2}{\ep_\la}+r^2|\ln r|\right),
\end{equation}
where $C>0$ is a constant independent of $\ep_\la$.

In particular, recalling \eqref{relpara}, we have the crucial estimate
\begin{equation}\label{key}
\la\mu_\la^2\approx \frac{C}{\ep_\la^2}\exp\left(\frac{A_1}{\varepsilon_\la} \left(\frac{1}{\xi (1)} -1 \right) \right) \rightarrow 0\quad \mbox{as}\ \la \rightarrow 0 \ (\mbox{and thus }\ep_\la\to0),
\end{equation}
where $A_1=\lim_{\la \to 0}A_{1,\ep_\la}>0$, $C>0$ is a constant independent of $\ep_\la$, and $\xi$ is the function defined in Lemma \ref{lemma:appendix1}, which satisfies $\xi(1)>1$ due to the fact that it is an strictly increasing function with $\xi (0)=1$.

\subsection{Construction of \texorpdfstring{$u_0$}{u0}}
Recall, by Lemma \ref{singreen}, that there exists 
\begin{equation}\label{deftilder}
\tilde{r}\approx\sqrt \varepsilon_\la 
\end{equation}
such that $u_2^\prime (\tilde{r})=0$.
We define 
\begin{equation}\label{defu_0}
u_0=V_{\mu_\la}+H_{0,\mu_\la},
\end{equation}
where $V_{\mu_\la}$ is the function defined in \eqref{defU_0} and $H_{0,\mu_\la}$ is the solution, for $\mu=\mu_\la$, to
\begin{equation}
\label{defH_0}
\left\{
\begin{array}{rcll}-\Delta H_{0,\mu} +H_{0,\mu}&=&-V_\mu& \mbox{in}\ (0,\tilde{r}) \\
H_{0,\mu}^\prime (\tilde{r})&=&-V_\mu^\prime (\tilde{r}).\end{array}\right.
\end{equation}
We introduced the function $H_{0,\mu_\la}$ in order to get a better matching between $u_0$ and $u_2$. 
We choose $\delta$ such that 
\begin{equation}
\label{defdelta}
2\delta < \tilde{r} \text{ and }\delta \approx \sqrt{\varepsilon_\la}. 
\end{equation}
Arguing in a similar way to the proof of \cite{MR2209293}*{Lemma 2.1}, we obtain the following estimates.
\begin{lem}
\label{expHo}
For any $\alpha \in \left(0,\frac{1}{2}\right)$, we have, for $r\in (0,\tilde{r})$,
\begin{equation}\label{RegEst}
H_{0,\mu_\la} (r)= H_{\ep_\la}(r)-  \ln ( 8\mu_\la^2 ) +O\left((\mu_\la^2\lambda)^\alpha\right)
\end{equation}
$C^{0,\gamma}(B_{\tilde r})-$uniformly, for $\gamma\in[0,1)$, where $H_{\ep_\la}(r)$ is defined in \eqref{regreen}. Moreover, \eqref{RegEst} holds uniformly in $C^1(B_{2\delta}\backslash B_\delta)$.
Finally, for $r\in (0,\tilde{r})$,
\begin{equation}
\label{estH_0}
|H_{0,\mu_\la}(r)|\leq C\left(\frac{r^2}{\ep_\la}+r^2|\ln r|+(\la\mu_\la^2 )^\alpha \right),
\end{equation}
where $C>0$ is a constant independent of $\ep_\la$.
\end{lem}

\begin{proof}
Let us consider the function $z(r)= H_{0,\mu_\la}(r) - H_{\ep_\la}(r) + \ln (8 \mu_\la^2 )$, which satisfies
$$\left\{\begin{array}{rcll}-\Delta z +z &=&-  \ln \dfrac{1}{(\la\mu_\la^2 +r^2 )^2} + \ln \dfrac{1}{r^4}& \mbox{in}\ (0,\tilde{r})\\
 z^\prime (\tilde{r})&=&\dfrac{4\tilde{r}}{\la \mu_\la^2  +\tilde{r}^2}-\dfrac{4}{\tilde{r}}.\end{array}\right. $$
By recalling \eqref{relpara} and that $\tilde{r}\approx \sqrt{\varepsilon}$, we deduce that
$$z^\prime (\tilde{r})= -\dfrac{4\la \mu_\la^2}{\tilde{r} (\la\mu_\la^2  + \tilde{r}^2)}=O((\la\mu_\la^2)^{\alpha}),$$
for any $\alpha \in \left(0,1\right)$. We set $f=-  \ln \dfrac{1}{(\la\mu_\la^2 +r^2 )^2} + \ln \dfrac{1}{r^4}$ and let $p>2$. We have
$$\int_{B_{\tilde{r}}} |f|^p dx =\int_{B_{\tilde{r}} \backslash B_{\sqrt{\lambda}\mu_\la }} |f|^p dx + \int_{ B_{ \sqrt{\lambda}\mu_\la }} |f|^p dx. $$
It is easy to see that
$$ \int_{ B_{\sqrt{\lambda}\mu_\la }} |f|^p dx\leq C  \lambda \mu_\la^2 |\ln (\la\mu_\la^2) |^{p},$$
and, using the fact that $|f(r)|\leq \dfrac{C \sqrt{\lambda}\mu_\la}{r}$, one gets
$$\int_{B_{\tilde{r}} \backslash B_{ \sqrt{\lambda}\mu_\la}} |f|^p dx \leq C  \lambda^{p/2} \mu_\la^p (\sqrt{\lambda}\mu_\la )^{2-p}\leq C\la \mu_\la^2.$$
Using elliptic regularity theory (see Lemma \ref{appregell}), we deduce that
$$\|z \|_{C^{0,\gamma}(B_{\tilde r})}\leq C(\lambda\mu_\la^2)^\alpha$$
for any $\gamma \in (0,1)$ and $\alpha \in \left(0,\frac12\right)$. 

On the other hand, for any $q\geq 2$, since $\delta\approx \sqrt{\varepsilon_\la}$, we have
$$\int_{B_{2\delta}\backslash B_{\delta}} |f|^q dx \leq C(\lambda\mu_\la^2)^{q/2} \delta^{2-q}\leq C(\lambda\mu_\la^2)^{q/2}\varepsilon_\la^{\frac12(2-q)}\leq C(\lambda\mu_\la^2)^{\alpha q},$$
for any $\alpha\in \left(0,\frac12\right)$. We deduce that
$$
\|z\|_{C^1(B_{2\delta}\backslash B_{\delta})}\leq C (\la\mu_\la^2)^\alpha.
$$
Finally, \eqref{estH_0} is a direct consequence of \eqref{RegEst} and \eqref{expansionH}.
\end{proof}
Thanks to the previous lemma, we are able to show that $u_0$ and $u_2$ are very close in $C^1-$norm sense in the interval $[\delta , 2\delta]$. 
\begin{lem}
\label{estlinea}
For $\delta \leq r \leq 2\delta$, we have
$$|u_0(r)-u_2(r)|= O\left((\la\mu_\la^2)^\alpha\right)\quad \mathrm{and}\quad |u_0^\prime (r)-u_2^\prime (r)|= O\left((\la\mu_\la^2)^\alpha\right)$$
for any $\alpha \in \left(0,\frac12\right)$.
\end{lem}
\begin{proof}
The proof is a direct consequence of Lemma \ref{expHo}. Indeed, by definition, for $r\in [\delta ,2\delta]$, we have
$$u_0(r)=V_{\mu_\la}(r)+H_{0,\mu_\la}(r)= \ln \dfrac{8\mu_\la^2}{(\lambda\mu_\la^2  +r^2 )^2} + H_{\ep_\la}(r)-\ln (8 \mu_\la^2) +O((\la\mu_\la^2)^\alpha) $$
and
$$u_2 (r)=- 4\ln r+H_{\ep_\la}(r).$$
It follows that
\begin{align*}
u_0(r) -u_2 (r)&=  -2 \ln \left(1+\dfrac{\la\mu_\la^2 }{r^2}\right) + O\left((\la\mu_\la^2)^\alpha\right)\\
&= O\left((\la\mu_\la^2)^\alpha\right). 
\end{align*}
Arguing in a similar way, one shows that
$$u_0^\prime (r)- u_2^\prime (r)=O\left(\dfrac{\la\mu_\la^2 }{\delta^3}\right)+ O\left((\la\mu_\la^2)^\alpha\right)= O\left((\la\mu_\la^2)^\alpha\right). $$
\end{proof}

\medskip 
We will now look for a solution of \eqref{ks} of the form $U+\phi$, where $\phi$ is a correction term. Let us observe that $U+\phi$ is a solution to \eqref{ks} if and only if $\phi$ solves  
\begin{equation}\label{eqdephi}
\left\{
\begin{array}{rcll}
L(\phi)&=&N(\phi)-R(U)&\mathrm{in}\ (0,1)\\
\phi'(0)&=&0&\\
\phi'(1)&=&0,&
\end{array}\right.
\end{equation}
where 
\begin{align}
L(\phi)&=-\Delta \phi +\phi -\lambda e^U \phi \label{defL}\\
N(\phi)&=\lambda ( e^{U+\phi}- e^U-e^U \phi) \label{defN}\\
R(U)&=-\Delta U +U -\lambda e^U\nonumber.
\end{align}

\section{The error estimate}\label{sec4}
In this section we estimate the terms $R(U)$ and $N(\phi)$. In order to take benefit of the estimates in \cite{PistoiaVaira2015}, we are going to work with the norm $\|\cdot\|_\ast$ (see \eqref{normstar}) which is,roughly speaking, a weighted $L^\infty-$norm on $B_{\frac12}$ and a $L^1-$norm elsewhere. We begin by estimating $N(\phi)$.

\begin{lem}
\label{estN}
There exists $C>0$ such that, for any $\beta>0$,
\begin{equation}\label{eN}
|N(\phi)|\leq C |\phi |^2 \left\{\begin{array}{cl}  \dfrac{1}{\la\mu_\la^2\left(1  +\left(\dfrac{r}{\sqrt{\lambda}\mu_\la}\right)^2\right)^2}& \mbox{if}\ r\leq \delta \\ \varepsilon_\la^\beta & \mbox{if}\ \delta \leq r \leq 1-2\delta_1 
 \end{array}\right.
\end{equation}
and
\begin{equation}
\label{estNext}
\left\|N(\phi)\right\|_{L^1\left(B_1 \backslash B_\frac{1}{2}\right)}\leq C\varepsilon_\la^{-1}  \|\phi \|_{L^\infty\left(B_1 \backslash B_\frac{1}{2}\right)}^2 .
\end{equation}
\end{lem}
\begin{proof}
First, using a Taylor's expansion, it is immediate to see that
$$|N(\phi)|\leq C \lambda e^U |\phi|^2.$$
Therefore, the proof reduces to estimate $e^U$. Let us observe that if $r\in [0,2\delta]$, using \eqref{estH_0} and \eqref{defdelta}, we have that
$$
\lambda e^{u_0}=\lambda e^{V_{\mu_\la}+H_{0,\mu_\la}} \leq C\dfrac{\mu_\la^2}{(\la\mu_\la^2 +r^2 )^2} \exp\left( \dfrac{r^2}{\varepsilon_\la }\right)\leq \frac{C}{\la\mu_\la^2\left(1  +\left(\dfrac{r}{\sqrt{\lambda}\mu_\la}\right)^2\right)^2}.$$

Besides, by definition of $u_2$, we know that it is decreasing in $(0,\tilde{r})$ and increasing elsewhere. Then, for $r\in [\delta, 1- 2\delta_1]$, we have
$$e^{u_2 (r)}\leq e^{u_2 (\delta)} +e^{u_2 (1-2\delta_1)}. $$
Making a Taylor's expansion, we obtain, for some $\theta \in (1-2\delta_1 ,1)$,
$$u_2 (1-2\delta_1)= u_2 (1)-2\delta_1 u_2^\prime (1)+2 \delta_1^2 u_2^{\prime \prime}(\theta) \leq \dfrac{\sqrt{2}}{\varepsilon_\la G_{\ep_\la}^\prime (1)}- \delta_1 u_2^\prime (1). $$
Thus, noting that $u_2'(1)\approx \sqrt2\ep_\la^{-1}$ and recalling \eqref{relpara} and \eqref{defdelta1}, we deduce that
$$\lambda e^{u_2 (1-2\delta_1)}\leq C \varepsilon_\la^{-2} e^{-\delta_1 u_2^\prime (1)  }\leq C \varepsilon_\la^\beta$$
for any $\beta>0$. On the other hand, recalling \eqref{defu_2}, we see that $e^{u_2 (\delta)} \leq \dfrac{C}{\delta^4}\leq C \varepsilon_\la^{-8}$. 
By noticing that $\lambda \varepsilon_\la^{-8}\leq C \varepsilon_\la^\beta $ for any $\beta>0$, we conclude that
$$
\lambda e^{u_2(r)}\leq C\varepsilon_\la^\beta
$$
for any $r\in [\delta,1-2\delta_1]$ and any $\beta>0$. Finally, by observing that for $r\in (\delta ,2\delta)$, 
$$\lambda e^U \leq \lambda \max \{e^{u_0},e^{u_2} \}\leq C\varepsilon_\la^\beta,$$
we deduce \eqref{eN}.

We refer to \cite{PistoiaVaira2015}*{Lemma  $4.3$} for the proof of \eqref{estNext}.

\end{proof}

Next, we estimate $R(U)$.
\begin{lem}
\label{estR}
Let $\alpha \in \left(0,\frac12\right)$. There exists $C>0$ such that
$$|R(U)| \leq C\left\{\begin{array}{cl} \dfrac{(\la\mu_\la^2)^\alpha+\dfrac{r^2}{\ep_\la}+r^2|\ln r|}{\la\mu_\la^2\left(1 +\left(\dfrac{r}{\sqrt{\lambda}\mu_\la}\right)^2 \right)^2}  & \mbox{if}\ r\leq \delta \\ \varepsilon_\la^\beta & \mbox{if}\ \delta \leq r \leq 1-2\delta_1 
\end{array}\right. $$
for any $\beta>0$, and
$$\left\|R(U)\right\|_{L^1\left( B_1 \backslash B_{\frac{1}{2}}\right)}\leq C\varepsilon_\la^{1+\sigma}$$
for some $\sigma>0$.
\end{lem}
\begin{proof}
First, we consider the case $r\leq \delta$ so that $U(r)=u_0(r) (r)=V_{\mu_\la}(r)+H_{0,\mu_\la}(r)$. Combining \eqref{defU_0}, \eqref{defH_0}, and \eqref{estH_0}, and making a Taylor's expansion, we infer that
\begin{align}
\label{juinestu0}
|R(u_0)|&=\left|-\Delta (V_{\mu_\la}+H_{0,\mu_\la})+V_{\mu_\la}+H_{0,\mu_\la} -\lambda e^{V_{\mu_\la}+H_{0,\mu_\la}}\right| \nonumber\\
&=\left|\lambda e^{V_{\mu_\la}} \left(1-e^{H_{0,\mu_\la}}\right)\right|\nonumber\\
&\leq C\la e^{V_{\mu_\la}}\left|H_{0,\mu_\la}\right|\nonumber\\
&\leq C\dfrac{(\la\mu_\la^2)^\alpha+\frac{r^2}{\ep_\la}+r^2|\ln r|}{\la\mu_\la^2\left(1  +\left(\dfrac{r}{\sqrt{\lambda}\mu_\la}\right)^2 \right)^2}.  
\end{align}

\noindent Next, when $2\delta \leq r\leq 1-2\delta_1$, we have $U(r)=u_2(r)$. Arguing as in the previous lemma, we obtain
\begin{equation}
\label{juinestu2}
|R(u_2(r))|=\lambda e^{u_2 (r)}\leq C\varepsilon_\la^\beta
\end{equation}
for any $\beta>0$.

On the other hand, it is proven in \cite{PistoiaVaira2015}*{Lemma $4.2$} that
\begin{equation}
\label{estborde11}
\| R(u_4)\|_{L^1\left(B_1 \backslash B_{1-\delta_1}\right)}\leq C \varepsilon_\la^{1+\sigma}\quad \mbox{for some } \sigma>0.
\end{equation}

Finally, we consider the two intermediate regimes. First, let us consider the case $\delta \leq r \leq 2\delta$. In this interval, we have $U(r)=u_1(r)$. Using \eqref{defu_1},  we get
\begin{align}
\label{juinestu1}
|R(u_1|)&=\left|\chi_1 R(u_0)+(1-\chi_1) R(u_2)-2 \chi_1^\prime (u_0^\prime - u_2^\prime)+(-\Delta \chi_1 +\chi_1) (u_0 -u_2)\right.\nonumber\\
&\left.\ \ \ \ + \lambda \chi_1 e^{u_0}+\lambda (1-\chi_1) e^{u_2}-\lambda e^{\chi_1 u_0 +(1-\chi_1) u_2}\right|\nonumber\\
&\leq |R(u_0)|+|R(u_2)|+C \left( \dfrac{|u_0^\prime - u_2^\prime|}{\delta}+\dfrac{|u_0 - u_2|}{\delta^2}\right)\nonumber\\
&\ \ \ \ + \lambda e^{u_2}+\left|\lambda e^{u_0}\left( e^{(\chi_1 -1) (u_0-u_2)}  -1 \right)\right|.
\end{align}
Using a Taylor's expansion and Lemma \ref{estlinea}, we have 
$$\left|\lambda e^{u_0}\left( e^{(\chi_1 -1) (u_0-u_2)}  -1 \right)\right| \leq \lambda e^{u_0}|u_0 -u_2| \leq  \la e^{u_0} (\la \mu_\la^2 )^\alpha. $$
Using Lemma \ref{estlinea} once again, we get
$$\dfrac{|u_0^\prime - u_2^\prime|}{\delta}+\dfrac{|u_0 - u_2|}{\delta^2}\leq C (\la \mu_\la^2)^\alpha \delta^{-2}.$$
Plugging these two last estimates into \eqref{juinestu1} and using \eqref{juinestu0} and \eqref{juinestu2}, we obtain
\begin{align*}
%\\
|R(u_1)|&\leq C\left( \sup_{\delta \leq r\leq 2\delta}
\dfrac{1}{\la \mu_\la^2 \left(1  +\left(\dfrac{r}{\sqrt{\lambda_\la}\mu_\la}\right)^2 \right)^2} +\varepsilon_\la^\beta +\dfrac{(\lambda \mu_\la^2)^\alpha}{\delta^2}\right)\\
&\leq C\left( \dfrac{\lambda \mu_\la^2} {\delta^4}+\varepsilon_\la^\beta +\dfrac{(\la \mu_\la^2)^\alpha}{\delta^2}\right) \leq C \varepsilon_\la^\beta 
\end{align*}

Finally, when $1-2\delta_1 \leq r \leq 1-\delta_1$, arguing as above, we have
\begin{align}
\label{juinestu3}
|R(u_3)|&=\left| \chi_3 R(u_2)+(1-\chi_3) R(u_4)-2 \chi_3^\prime (u_4^\prime - u_2^\prime)+(-\Delta \chi_3 +\chi_3) (u_4 -u_2)\right.\nonumber\\
&\ \ \ \ \left.+ \lambda \chi_3 e^{u_4}+\lambda (1-\chi_3) e^{u_2}-\lambda e^{\chi_3 u_4 +(1-\chi_3) u_2}\right|\nonumber\\
&\leq |R(u_2)|+|R(u_4)|+C \left( \dfrac{|u_4^\prime - u_2^\prime|}{\delta_1}+\dfrac{|u_4 - u_2|}{\delta_1^2}\right)\nonumber\\
&\ \ \ \ + \lambda e^{u_2}+\lambda e^{u_4}\left| u_4 - u_2  \right|.
\end{align}
Using Lemma \ref{lemrelu2u4} and the definition of $\delta_1$ given in \eqref{defdelta1}, we obtain
$$\int_{1-2\delta_1}^{1-\delta_1} \left( \dfrac{|u_4^\prime - u_2^\prime|}{\delta_1}+\dfrac{|u_4 - u_2|}{\delta^2_1}\right) rdr =O(\varepsilon_\la^{1+\sigma})$$
and
$$\int_{1-2\delta_1}^{1-\delta_1}\lambda e^{u_4}\left| u_4 - u_2  \right| r dr=O(\varepsilon_\la^{1+\sigma}).  $$

Thanks to \eqref{juinestu2} and \eqref{estborde11}, we see that
$$\int_{1-2\delta_1}^{1-\delta_1}  \left( R(u_2)+R(u_4) +\lambda e^{u_2} \right) r dr =O(\varepsilon_\la^{1+\sigma}).$$
Plugging the three previous estimates into \eqref{juinestu3}, we obtain
$$\| R(u_3)\|_{L^1 \left( B_1 \backslash B_{\frac{1}{2}}\right) }\leq C\varepsilon_\la^{1+\sigma}.$$
This concludes the proof of the lemma.
\end{proof}

\section{Invertibility of the linearized operator}\label{sec5}
In this section we develop an invertibility theory for the operator $L$ defined in \eqref{defL}. To do so, we use ideas from \cites{delPinoKowalczykMusso,MR2209293,delPinoRoman,PistoiaVaira2015}. First, we define the norms
\begin{align}
\label{normstar}
\|u \|_\ast =&\max \left\{|\ln (\lambda\mu_\la^2) |\|\tilde{\chi}_1 u \|_{\star}, \|\tilde{\chi}_2 u \|_{L^1(B_1\setminus B_{1/4})} \right\}
\end{align}
and
\begin{align*}
\|u \|_{\ast\ast} =&\max \left\{\|\tilde{\chi}_1 u \|_{\star}, \|\tilde{\chi}_2 u \|_{L^1(B_1\setminus B_{1/4})} \right\},
\end{align*}
where 
$$\tilde{\chi}_1 (r)= \left\{\begin{array}{ll}1 &\mbox{if } r\leq \frac12\\[.2cm]  0& \mbox{if } r\geq \frac34 \end{array} \right., \quad 
  \tilde{\chi}_2 (r)= \left\{\begin{array}{ll}1 &\mbox{if } r\geq \frac12\\[.2cm]  0& \mbox{if } r\leq \frac14 \end{array} \right.,$$
and
$$\| u \|_{\star} = \sup_{r\in\left[0,\frac34\right]} \frac{\lambda\mu_\la^2 |u(r)|}{\lambda\mu_\la^2+ \left(1+\frac{r}{\sqrt{\lambda}\mu_\la}\right)^{-2- \nu }}=\sup f_\lambda (r) |u(r)|$$
for some $\nu \in (0,1)$. 
 
The following proposition is the main result of this section. 
\begin{prop}\label{teolineal}
There exist positive constants $\lambda_0$ and $C$ such that for any $\lambda \in (0,\lambda_0)$ and for any $h\in L^\infty(B_1)$, there exists a unique radial function $\phi \in W^{2,2}(B_1)$ solution to
\begin{equation}\label{lineari0}\left\{\begin{array}{rcll} 
L(\phi )&=&h&\mathrm{in}\ B_1\\
\phi^\prime(1)&=&0,&
\end{array}\right. \end{equation}
which satisfies 
\begin{equation}\label{estimateLinear}
\left\|\phi\right\|_{L^\infty(B_1)}\leq C \left\|h\right\|_{\ast}.
\end{equation}
\end{prop}

Rather than directly proving this statement, we first provide a priori estimates for the solution to \eqref{lineari0} when $\phi$ is orthogonal to
$$
Z_0(r)=\frac{r^2-\lambda\mu_\la^2}{r^2+\lambda\mu_\la^2}.
$$
This function is a solution to
\begin{equation}\label{eqDilations}
-\Delta Z_0 =\frac{8\la\mu_\la^2}{(\lambda\mu_\la^2+r^2)^2}Z_0,
\end{equation}
that is, the linearization of the equation $-\Delta v=e^v$ around the radial solution 
$$
v(r)=V_{\mu_\la} (r) +\ln \lambda =\ln \dfrac{8\lambda\mu_\la^2}{(\lambda\mu_\la^2+|r|^2)^2}.
$$
It is well-known that the only bounded radial solutions to \eqref{eqDilations} are multiples of $z_0$ (see \cite{MR1885666}*{Lemma 2.1}). 

\medskip
Consider a large but fixed number $R_0>0$ and a radial smooth cut-off function $\chi_\la(r)$ such that $\chi_\la(r)=1$ if $r\leq R_0\sqrt \lambda\mu_\la$ and $\chi_\la(r)=0$ if $r>(R_0+1)\sqrt\lambda\mu_\la$. 

\begin{lem}\label{teolineal2}
There exist positive constants $\lambda_0$ and $C$ such that for, any $\lambda \in (0,\lambda_0)$, the unique radial solution $\phi \in W^{2,2}(B_1)$ to
\begin{equation}\label{lineari}\left\{\begin{array}{rcll} 
L(\phi )&=&h& \mathrm{in}\ B_1\\ 
\phi^\prime (1) &=&0&\\
\displaystyle \int_{B_1} \chi_\la Z_0 \phi \ dx&=&0 \end{array}\right. \end{equation}
satisfies
$$\left\|\phi\right\|_{L^\infty(B_1)}\leq C \left\|h\right\|_{\ast\ast}.$$
\end{lem}

\begin{proof}
Assume towards a contradiction that there exist a sequence of positive numbers $\lambda_n \rightarrow 0$ and a sequence of solutions $\phi_n$ to \eqref{lineari} such that 
\begin{equation}
\label{contre1}
\left\|\phi_n\right\|_{L^\infty (B_1)}=1\quad \mathrm{and}\quad \left\|h_n\right\|_{\ast\ast} \underset{n\rightarrow \infty}{\longrightarrow} 0.
\end{equation}
We denote by $\varepsilon_n$ (resp. $\mu_n$) the sequence defined by \eqref{relpara} (resp. \eqref{defmu}) with $\la=\la_n$. Also, we let $U_n$ denote the first approximation for $\la=\la_n$ and $u_0^n$ (resp. $u_2^n$) the sequence defined by \eqref{defu_0} (resp. \eqref{defu_2}) with $\la=\la_n$.

Our goal is to prove that $\phi_n (r)=o_n(1)$, for any $r\in [0,1]$, which yields to a contradiction with \eqref{contre1}, where here and in the rest of the proof $o_n(1)$ denotes a function $f_n(r)$ such that $\displaystyle\lim_{n\rightarrow \infty} f_n(r)=0$ uniformly in $r$. We split the proof in $4$ steps.

\bigskip
\noindent{\bf Step 1.} There holds $\phi_n (r)=o_n (1)$ on compact subsets of $(0,1)$.

\smallskip
First, it is easy to see that $\la_ne^{U_n}=o_n(1)$ on compact subsets of $(0,1)$. Since, by assumption, $\left\|h_n\right\|_{**} \rightarrow 0$, we deduce that, up to subsequence, $\phi_n\to \hat\phi$ $C_2-$uniformly on compact subsets of $(0,1)$, where $\hat \phi$ is a radial bounded solution to 
\begin{equation}
\label{limitint}
\left\lbrace\begin{array}{rcll}
-\Delta \hat \phi +\hat \phi &=&0& \mbox{in } B_1\setminus \{0\}\\
 \hat \phi' (1)&= & 0.
\end{array}\right.
\end{equation}
We claim that $\hat{\phi}\equiv 0$, which in turn implies that $\phi_n (r)=o_n (1)$ on compact subsets of $(0,1)$. To prove the claim, let us observe that equation \eqref{limitint} corresponds to the modified zero-order Bessel differential equation, whose solution is given by
$$
\hat \phi(r)= C_1 \xi(r)+C_2 \zeta(r),
$$
where $C_1,C_2$ are constants, $\xi(r)$ is the zero-order modified Bessel function of the first kind, and $\zeta(r)$ is the zero-order modified Bessel function of the second kind defined in Lemma \ref{lemma:appendix1}. Since $\zeta(r)$ becomes unbounded as $r\to 0$ and $\hat \phi$ is bounded, we deduce that $C_2=0$. Moreover, since $\xi'(1)\neq 0=\hat\phi'(1)$, we have $C_1=0$. Hence, $\hat \phi=0$.

\bigskip
\noindent{\bf Step 2.} We have that $\phi_n (r)=o_n (1)$ for $r$ close to $1$.

\smallskip 
We set $\psi_n (s)= \phi_n (\varepsilon_n s+1)$ for $s\in[-\varepsilon_n^{-1},0]$. Then, since $\psi_n$ is bounded, by arguing as in \cite{PistoiaVaira2015}*{Proposition 5.1} it is possible to show that $\psi_n \rightarrow \psi$ $C^2$--uniformly on compact subsets of $(-\infty , 0]$ where $\psi$ satisfies
$$
\left\{
\begin{array}{rcll}
-\psi^{\prime \prime}&=&e^{\psi}&\mbox{in}\ \R^{-}\\ \psi^\prime (0)&=&0\\
\|\psi\|_{L^\infty}&\leq &1.\end{array}\right.$$
We know (see \cite{Grossi2006}) that any solution $\psi$ to $-\psi^{\prime \prime}= e^{\psi}$ is of the form
$$\psi (s)=a\dfrac{e^{\sqrt{2} s}-1}{e^{\sqrt{2}s}+1}+b \left(-2+\sqrt{2} s \dfrac{e^{\sqrt{2} s}-1}{e^{\sqrt{2}s}+1}\right)$$
for some $a,b \in \R$. Since $ \|\psi\|_{\infty}\leq 1$, we have $b=0$, and since $\psi'(0)=0$, we get $a=0$. Hence, $\psi\equiv 0$.

\medskip
Let us now consider the radial Green's function $G(r,t)$ associated to the operator 
$$(-\Delta \cdot +\ \cdot)
$$ 
satisfying $G\left(r,\frac12\right)=G'(r,1)=0$ and singular at the point $r\in \left(\frac12,1\right)$. Using Green's formula, we have that, for $\frac12 <r< 1$,
\begin{align*}
\phi_n (r)-G^\prime \left(r,\frac12 \right)\phi_n \left(\frac12\right)&= \int_{\frac12}^1 G(r,t) h_n (t) dt+ \lambda_n\int_{\frac12}^1 G(r,t)  e^{U_n}\phi_n (t)dt\\
&=\int_{\frac12}^1 G(r,t) h_n (t) dt+ G(r,1) \varepsilon_n \lambda_n \int^0_{-\frac1{2\varepsilon_n}} e^{U_n (\varepsilon_n s+1)} \psi_n (s) ds\\
&+   \varepsilon_n \lambda_n \int^0_{-\frac1{2\varepsilon_n}} (G(r,\varepsilon_n s +1)-G(r,1)) e^{U_n (\varepsilon_n s+1)} \psi_n (s) ds.
\end{align*}
From {\bf Step 1} we know that $\phi_n\left(\frac12\right)=o_n(1)$. Combining this with the fact that $G$ is bounded and $\|h_n\|_{**}\to 0$ as $n \to \infty$, we have
$$
G^\prime \left(r,\frac12 \right)\phi_n \left(\frac12\right)+\int_{\frac12}^1 G(r,t) h_n (t) dt=o_n (1).
$$ 
Arguing as in \cite{PistoiaVaira2015}*{Proposition 5.1}, one shows that
$$ 
\varepsilon_n \lambda_n \int^0_{-\frac1{2\varepsilon_n}} (G(r,\varepsilon_n s +1)-G(r,1)) e^{U_n (\varepsilon_n s+1)} \psi_n (s) ds=o_n (1).
$$ 
Hence,
$$
\phi_n(r)=C_nG(r,1)+o_n(1),
$$
where $C_n=\varepsilon_n \lambda_n \int^0_{-\frac1{2\varepsilon_n}} e^{U_n(\varepsilon_n s+1)} \psi_n (s) ds$. 

Since $\phi_n(1)=\psi_n(1)=o_n(1)$ and $\lim_{r\to 1}G(r,1)\neq 0$, we deduce that $C_n=0$. Hence, $\phi_n=o_n(1)$ for $r\in \left[\frac12,1\right]$. 

\bigskip
In the following steps it is convenient to work with rescaled variables. We set 
$$
s=\frac{r}{\sqrt \lambda_n\mu_n}\quad \mbox{for }r\in[0,1],
$$
and define 
\begin{align*}
\tilde \phi_n (s)&=\phi_n(\sqrt{\lambda_n}\mu_n s),\\
\tilde{U}_n(s)&=U_n (\sqrt{\lambda_n}\mu_n s) + 2\ln (\lambda_n \mu_n),\\
\tilde h_n(s)&=\lambda_n\mu_n^2h_n(\sqrt{\lambda_n}\mu_n s).
\end{align*}
Letting $\tilde L(\cdot)=-\Delta \cdot +\lambda_n\mu_n^2 \cdot -e^{\tilde U_n}\cdot $, it easy to see that $\tilde \phi_n$ satisfies
\begin{equation}\label{eqphi} 
\tilde{L} (\tilde\phi_n(s))=\tilde h_n(s).
\end{equation}
We also define (with some abuse of notation)
\begin{equation}\label{scalednormstar}
\|\tilde h\|_{\star}:=\sup_{s\in \left[0,\frac1{2\sqrt \lambda_n \mu_n }\right]}\frac{\tilde h(s)}{\lambda_n\mu_n^2+(1+s)^{-2-\nu}}=\|h1_{\{r\leq \frac12\}}\|_{\star},
\end{equation}
for functions $\tilde h$ defined in the rescaled variable. 

\bigskip
{\bf Step 3.} Up to subsequence, we have that $\tilde \phi_n\to 0$ as $n\to\infty$ uniformly over compact sets of $\R^2$.

\smallskip
Recalling \eqref{estH_0}, elliptic estimates applied to \eqref{eqphi} imply that, up to subsequence, $\tilde{\phi}_n$ converges uniformly over compact sets of $\R^2$ to a bounded solution $\tilde{\phi}$ to
$$
-\Delta \tilde{\phi}= \frac{8}{(1+s^2)^2} \tilde{\phi}\quad \mbox{in } \R^2.
$$
By \cite{MR1885666}*{Lemma 2.1}, we deduce that there exists a constant $C_0$ such that $\tilde \phi= C_0 \tilde Z_0(s)$, where 
$$\tilde Z_0(s)=\frac{s^2-1}{s^2+1}=Z_{0,n}(\sqrt{\lambda_n}\mu_n s)\quad\mathrm{with}\ 
Z_{0,n}(r)=\frac{r^2-\lambda_n\mu_n^2}{r^2+\lambda_n\mu_n^2}.
$$ 
Let $\tilde \chi(s)=\chi_{\la_n}(\sqrt{\la_n}\mu_n s)$ where $\chi_{\la}$ is defined just before the statement of the lemma. Notice that it does not depend on $n$. The orthogonality condition satisfied by $\phi_n$ yields
$$0=\int_{B_1} \chi_{\la_n} Z_{0,n} \phi_n dx=\lambda_n\mu_n^2\int_{B_{1/(\sqrt{\lambda_n}\mu_n)}} \tilde \chi \tilde Z_0 \tilde \phi_n dx.$$ 
Passing to the limit $n\to\infty$, we find
$$
\int_{\R^2} \tilde \chi \tilde Z_0 \tilde \phi dx=0,
$$
which implies $C_0=0$. The result thus follows. 

\bigskip
{\bf{Step 4.}} We have that $\phi_n (r)=o_n(1)$ for $r$ close to $0$.

\smallskip
This is based on a maximum principle argument. Let us show that there exists a constant $C>0$, independent of $n$, such that
\begin{equation}\label{InnerNormEstimate}
\|\tilde{\phi}_n\|_{L^\infty\left(B_{\tau/(\sqrt{\lambda_n}\mu_n)}\right)}\leq C\left[\sup_{s\leq R} |\tilde{\phi}_n(s)| +\|\tilde{h}_n \|_{\star} +|\phi_n(\tau)|\right],
\end{equation}
where $R>0$ is a large but fixed number and $\tau$ is a small but fixed number. 

To prove this, we need the following version of the maximum principle. We claim that 
there exists a fixed number $R_1>0$ such that for all $R>R_1$, 
$$\mbox{if }\tilde L(\varphi)>0 \mbox{ in } A_n:=B_{\tau/(\sqrt{\la_n}\mu_n)}\setminus B_R \mbox{ and }\varphi\geq 0 \mbox{ on } \partial A_n, \mbox{ then }\varphi\geq 0 \mbox{ in }A_n.
$$
To prove this, we consider the function $\varphi_0(s)=1-\frac{1}{s^{\nu}}$. Observe that 
$$
-\Delta \varphi_0= \nu(\nu+2) \frac1{s^{2+\nu}}.
$$
Letting $\tilde u_0^n(s)=u_0^n(\sqrt{\la_n}\mu_ns)+2\ln(\la_n\mu_n)$ and $\delta_n$ be defined by \eqref{defdelta} with $\la=\la_n$, and using \eqref{estH_0}, we deduce that, for any $s\in[0,2\delta_n/\sqrt{\la_n}\mu_n]$,
\begin{equation}\label{es0}
e^{\tilde u_0^n(s)}=\frac{8}{(1+s^2)^2}e^{H_{0,\mu_{\la_n}}(\sqrt{\la_n}\mu_ns)}\leq \frac{C}{(1+s^2)^2}.
\end{equation}
On the other hand, letting $\tilde u_2^n(s)=u_2^n(\sqrt{\la_n}\mu_ns)+2\ln(\la_n\mu_n)$,
and $\tilde r_n$ be defined by \eqref{deftilder} with $\la=\la_n$, and recalling that $u_2^n$ is increasing for $r>\tilde r_n$ and \eqref{defmu},
we have, for any $s\in[\tilde r_n /(\sqrt{\la_n}\mu_n),\tau/(\sqrt{\la_n}\mu_n)]$,
\begin{equation}\label{es1}
e^{\tilde u_2^n(s)}\leq e^{\tilde u_2^n (\tau/(\sqrt{\la_n}\mu_n))}=8 \left(\frac{\sqrt{\la_n}\mu_n}{\tau}\right)^4\exp\left( 
H_{\ep_{\la_n}}(\tau)-H_{\ep_{\la_n}}(0)\right).
\end{equation} 
Choosing $\tau$ small but fixed so that \eqref{expansionH} holds, we have that 
$$
\exp\left( 
H_{\ep_{\la_n}}(\tau)-H_{\ep_{\la_n}}(0)\right)\leq C\exp\left(\frac{\tau^2}{\ep_\la}\right).
$$
In view of \eqref{key}, choosing $\tau$ smaller if necessary, we have that
$$
\frac{\sqrt{\la_n}\mu_n}{\tau}\exp\left( 
H_{\ep_{\la_n}}(\tau)-H_{\ep_{\la_n}}(0)\right)\leq C \frac{\sqrt{\la_n}\mu_n}{\tau}\exp\left(\frac{\tau^2}{\ep_\la}\right)\leq C.
$$
Plugging this into \eqref{es1}, we deduce that 
\begin{equation}\label{es2}
e^{\tilde u_2^n(s)}\leq e^{\tilde u_2^n (\tau/(\sqrt{\la_n}\mu_n))}\leq C\left(\frac{\sqrt{\la_n}\mu_n}{\tau}\right)^3\leq \frac{C}{s^3}.
\end{equation}
Besides, arguing similarly, it is easy to check that, for any $s\in[\delta_n/\sqrt{\la_n}\mu_n,\tilde r_n/\sqrt{\la_n}\mu_n]$,
$$
e^{\tilde u_2^n(s)}\leq \frac{C}{s^4}.
$$
From this, \eqref{es2}, and \eqref{es0}, we deduce that
$$
e^{\tilde U_n(s)}\leq \frac{C}{s^3}.
$$

Therefore, for any $R$ is sufficiently large, we have
$$
\tilde L(\varphi_0)=-\Delta \varphi_0 +\la_n\mu_n^2\varphi_0-e^{\tilde U_n}\varphi_0\geq \frac{\nu^2}2\frac{1}{s^{2+\nu}}+\frac12 \la_n\mu_n^2>0\quad \mbox{in } A_n
$$
and $\varphi_0>0$ on $\partial A_n$. The claim thus follows. 

\medskip
Let us now prove \eqref{InnerNormEstimate}. We define 
$$
\bar{\phi}_n=C_1\left[\max_{s\in (0,R)} |\tilde{\phi}_n (s)|+\|\tilde{h}_n\|_{\star}+|\phi_n(\tau)| \right]\varphi_0
$$
for a constant $C_1$ independent of $n$ and $R>R_1$. Observe that if $C_1\geq\dfrac{4}{\nu^2}$ then
$$
\tilde L(\bar \phi_n)\geq 2\|\tilde{h}_n\|_{\star} (s^{-2-\nu}+\lambda_n\mu_n^2)\geq |\tilde h_n|\frac{2(s^{-2-\nu}+\lambda_n\mu_n^2)}{\left((1+s)^{-2-\nu}+\lambda_n\mu_n^2\right)}\geq |\tilde h_n|=|\tilde L(\tilde \phi_n)|
$$ 
in $A_n$, since $\dfrac{2(s^{-2-\nu}+\lambda_n\mu_n^2)}{\left((1+s)^{-2-\nu}+\lambda_n\mu_n^2\right)}\geq 1$ for $s\in[R,+\infty)$ (taking $R$ larger if necessary).
On the other hand, if $C_1\geq (1-R^{-\nu})^{-1}$, we have
$$
\bar \phi_n \geq |\tilde \phi_n|\quad \mbox{on }\partial A_n.
$$
Applying the maximum principle and observing that $|\varphi_0|\leq 1$, we are led to \eqref{InnerNormEstimate}.

Finally, by noting that $\|\tilde h_n\|_\star\leq \|h_n\|_{**}=o_n(1)$ (by \eqref{contre1}), $\phi_n(\tau)=o_n(1)$ (by {\bf Step 1}), and $\max_{s\in (0,R)}|\hat\phi_n(s)|=o_n(1)$ (by {\bf Step 3}), we conclude that $\|\phi_n\|_{L^\infty(B_\tau)}=o_n(1)$. From this and {\bf Step 1 and Step 2}, we deduce that $\|\phi_n\|_{L^\infty(B_1)}=o_n(1)$ which contradicts the fact that $\|\phi_n\|_{L^\infty(B_1)} =1$. This completes the proof.
\end{proof}
We are now ready to prove Proposition \ref{teolineal}.

\begin{proof}[Proof of Proposition \ref{teolineal}]
We reuse the notation introduced in the proof of the previous lemma. For a scaled function $\tilde g (s)=\lambda\mu_\la^2 g(\sqrt \lambda\mu_\la s)$ with $s=r/(\sqrt \lambda\mu_\la )$, we define
\begin{equation}\label{scalednormast}
\|\tilde g\|_{**}:= \|g\|_{**}.
\end{equation}
Let $R>R_0+1$ be a large fixed number, $\tau>0$ as in {\bf Step 4} of the proof of Lemma \ref{teolineal2}, and $\hat z_0$ be the solution to
$$
\left\{
\begin{array}{rcll}
-\Delta \hat z_0 &=& \dfrac{8}{(1+s^2)^2}\hat z_0&\mbox{in } B_{\tau/(\sqrt\lambda\mu_\la) }\backslash B_R\\
\hat z_0(R)&=&\tilde Z_0(R)&\\
\hat z_0 \left(\tau/(\sqrt\lambda\mu_\la)\right)&=&0,&
\end{array}\right.
$$
where $\tilde{Z}_0$ is the function defined in {\bf Step 3} of the proof of Lemma \ref{teolineal2}. A direct computation shows that 
$$
\hat z_0(s)=\tilde Z_0(s)\left[1- \dfrac{\displaystyle \int_R^s \dfrac{dt}{t\tilde Z_0^2(t)} }{\displaystyle\int_R^{\tau/(\sqrt \lambda \mu_\la )} \dfrac{dt}{t\tilde Z_0^2(t)}}   \right].
$$
We also let $\hat z_1$ be the solution to 
$$
\left\{
\begin{array}{rcll}
-\Delta \hat z_1+\la\mu_\la^2\hat z_1 &=& \la\mu_\la^2\hat z_0&\mbox{in } B_{\tau/(\sqrt\lambda\mu_\la) }\backslash B_R\\
\hat z_1(R)&=&0&\\ 
\hat z_1 \left(\tau/(\sqrt\lambda\mu_\la)\right)&=&0.&
\end{array}\right.
$$
Elliptic estimates immediately yield that 
\begin{equation}\label{boundz1}
\|\hat z_1\|_{C^2\left(B_{\tau/(\sqrt\lambda\mu_\la) }\backslash B_R\right)}\leq C\la\mu_\la^2.
\end{equation}

We consider smooth cut-off functions $\eta_1(s)$ and $\eta_2(s)$ with the following properties: $\eta_1(s)=1$ for $s<R$, $\eta_1(s)=0$ for $s>R+1$, $|\eta'_1(s)|\leq 2$, $\eta_2(s)=1$ for $s<\tau/(2\sqrt\lambda\mu_\la )$, $\eta_2(s)=0$ for $s>\tau/(\sqrt\lambda\mu_\la )$, $|\eta'_2(s)|\leq C\sqrt \lambda\mu_\la$, and $|\eta''_2(s)|\leq C\lambda\mu_\la^2$.
We also define the test function 
$$
\tilde z_0=\eta_1 \tilde Z_0 +(1-\eta_1)\eta_2(\hat z_0-\hat z_1).
$$
Let $\phi$ be a solution to \eqref{lineari0}. As previously, we denote $\tilde \phi (s)= \phi (\sqrt{\lambda}\mu_\la s )$ and we let $\tilde{\chi}(s)=\chi_\la (\sqrt{\lambda}\mu_\la s)$ where $\chi_\la$ is defined before the statement of Lemma \ref{teolineal2}. Next, we modify $\tilde \phi$ so that the orthogonality condition with respect to $\tilde z_0$ is satisfied. We let
\begin{equation}\label{modif}
\hat \phi =\tilde \phi +A\tilde z_0,
\end{equation}
where the number $A$ is such that
$$
A\int_{B_{1/(\sqrt\lambda\mu_\la )}} \tilde{\chi} |\tilde{z}_0|^2 dx+\int_{B_{1/(\sqrt\lambda\mu_\la )}} \tilde{\chi} \tilde{z}_0 \tilde{\phi} dx =0.
$$
Then 
\begin{equation}\label{ModifiedProblem}
\tilde L(\hat \phi)=\tilde h+A\tilde L(\tilde z_0),
\end{equation}
and $\displaystyle\int_{B_{1/(\sqrt\lambda\mu_\la )}} \tilde{\chi} \tilde{z}_0 \hat \phi dx =0.$
Recalling \eqref{scalednormast}, Lemma \eqref{teolineal2} yields
\begin{equation}
\label{ModifiedEst}
\|\hat \phi\|_{L^\infty \left((B_{1/(\sqrt\lambda\mu_\la)}\right) }\leq C \left[\|\tilde h\|_{ \ast\ast}+|A|\|\tilde L(\tilde z_0)\|_{\ast\ast}    \right].
\end{equation}
Observe that $\tilde z_0=0$ for $s>\tau/(\sqrt\lambda\mu_\la )$. Thus, recalling \eqref{scalednormstar}, we have
$$
\|\tilde L(\tilde z_0)\|_{**}=\|\tilde L(\tilde z_0)\|_{\star}.
$$

Let us now estimate the size of $|A|$ and $\|\tilde L(\tilde z_0)\|_\star$.
Testing equation \eqref{ModifiedProblem} against $\tilde z_0$ and integrating by parts, we find
$$
\langle \hat \phi, \tilde L(\tilde z_0)\rangle = \langle \tilde h,\tilde z_0\rangle+A\langle \tilde L(\tilde z_0),\tilde z_0\rangle,
$$
where $\langle f,g \rangle =\displaystyle\int_{B_{1/(\sqrt\lambda\mu_\la)}}fg dx$. Combining this with \eqref{ModifiedEst},
$$
\int_{B_{1/(\sqrt\lambda\mu_\la )}} |\hat \phi||\tilde L( \tilde z_0) | dx \leq C\|\hat \phi \|_\infty \|\tilde L(\tilde z_0)\|_{\star},\  \text{ and} \quad\int_{B_{1/(\sqrt\lambda\mu_\la)}} |\tilde h| |\tilde z_0| dx\leq  C\|\tilde h\|_{\star},
$$
we are led to
\begin{equation}\label{estimateA}
A\langle \tilde L(\tilde z_0),\tilde z_0\rangle\leq C\|\tilde h\|_{\star}\left[1+\|\tilde L(\tilde z_0)\|_{\star} \right]+C|A|\|\tilde L(\tilde z_0)\|_{\star}^2.
\end{equation}
We next measure the size of $\|\tilde L(\tilde z_0)\|_{\star}$. We have 
\begin{align}
\label{traine1}
\tilde L(\tilde z_0)=&\eta_1\left(\frac{8}{(1+s^2)^2}\tilde Z_0 +\mu_\la^2 \lambda \tilde{Z}_0 -e^{\tilde U}\tilde Z_0\right)+ (1-\eta_1)\eta_2\left(\frac{8}{(1+s^2)^2}\tilde{z}_0-e^{\tilde U} (\hat z_0-\hat z_1) \right) \\
&+2\nabla \eta_1 \nabla (\eta_2 (\hat{z}_0-\hat z_1)  -\tilde Z_0)+ \Delta \eta_1 (\eta_2 (\hat{z}_0-\hat z_1)-\tilde Z_0 )\nonumber\\
&-2(1-\eta_1) \nabla \eta_2 \nabla (\hat z_0-\hat z_1)-(1-\eta_1)\Delta \eta_2 (\hat z_0-\hat z_1),\nonumber
\end{align}
where $\tilde U$ denotes, as in the proof of the previous lemma, the first approximation of the solution in the rescaled variable.
 
In the support of $\eta_1$, we have
$$
\frac{8}{(1+s^2)^2}\tilde Z_0-e^{\tilde U}\tilde Z_0=\frac{8}{(1+s^2)^2}\left(1-e^{\tilde H_{0,\mu_\la}}\right)\tilde Z_0,
$$
which combined with \eqref{estH_0} gives
$$
\left\|\eta_1\left(\frac{8}{(1+s^2)^2}\tilde Z_0 +\la \mu_\la^2 \tilde{Z}_0-e^{\tilde U}\tilde Z_0\right)\right\|_\star \leq C(\la\mu_\la^2)^\alpha.
$$
If $R\leq s\leq \delta/(\sqrt\la\mu_\la)$ (recall \eqref{defdelta}), we have that
$$
\frac{8}{(1+s^2)^2}\tilde{z}_0-e^{\tilde U} (\hat z_0-\hat z_1) =\frac{8}{(1+s^2)^2}\left(1-e^{\tilde H_{0,\mu_\la}}\right)\hat z_0+\frac{8}{(1+s^2)^2}e^{\tilde H_{0,\mu_\la}}\hat z_1.
$$
Therefore, using once again \eqref{estH_0} and \eqref{boundz1}, we deduce that
$$
\left\|\left(\frac{8}{(1+s^2)^2}\tilde{z}_0-e^{\tilde U} (\hat z_0-\hat z_1)\right)1_{\{R\leq s\leq \delta/(\sqrt\la\mu_\la)\}}\right\|_\star \leq C (\la\mu_\la^2)^\alpha.
$$
Arguing as in {\bf Step 4} of the proof of Lemma \eqref{teolineal2}, we deduce that, for any $s\in[\delta/(\sqrt\la \mu_\la),\tau/(\sqrt\la \mu_\la)]$,
$$
e^{\tilde U}\leq \frac{C}{s^3}.
$$
Hence,
$$
\left\|\left(\frac{8}{(1+s^2)^2}\tilde{z}_0-e^{\tilde U} (\hat z_0-\hat z_1)\right)1_{\{\delta/(\sqrt\la \mu_\la)\leq s\leq \tau/(\sqrt\la\mu_\la)\}}\right\|_\star \leq C (\la\mu_\la^2)^\alpha,
$$
and therefore
\begin{multline}\label{ee1}
\Bigg\|\eta_1\left(\frac{8}{(1+s^2)^2}\tilde Z_0+\la \mu_\la^2 \tilde{Z}_0-e^{\tilde U}\tilde Z_0\right)\\ + (1-\eta_1)\eta_2\left(\frac{8}{(1+s^2)^2}\tilde{z}_0-e^{\tilde U} (\hat z_0-\hat z_1) \right)\Bigg\|_\star \leq C(\la\mu_\la^2 )^\alpha.
\end{multline}

On the other hand, it is easy to see that, for $s\in (R,R+1)$, 
$$
|\tilde{Z}_0-\hat z_0|=\ \displaystyle \vline \tilde{Z}_0 \dfrac{\displaystyle \int_R^s \dfrac{dt}{t\tilde Z_0^2(t)} }{\displaystyle\int_R^{\tau/(\sqrt\la\mu_\la )}\dfrac{dt}{t\tilde Z_0^2(t)}}\vline \leq C|\ln (\lambda\mu_\la^2 )|^{-1}\quad \mathrm{and}\quad |\tilde{Z}_0'-\hat z_0'|\leq C|\ln (\lambda\mu_\la^2 )|^{-1}.
$$
Besides, for $s\in \left(\tau/(2\sqrt\lambda\mu_\la ),\tau/(\sqrt\lambda\mu_\la) \right)$, we have 
\begin{equation}
\label{traine2}
|\hat z_0|\leq C|\ln (\lambda\mu_\la^2 )|^{-1} \text{ and} \quad |\hat z_0'|\leq C \sqrt{\lambda}\mu_\la |\ln (\lambda\mu_\la^2 )|^{-1}.
\end{equation}
We then easily deduce, using \eqref{boundz1}, that 
\begin{align*}
&\|2\nabla \eta_1 \nabla (\eta_2 (\hat{z}_0-\hat z_1)  -\tilde Z_0)+ \Delta \eta_1 (\eta_2 (\hat{z}_0-\hat z_1)-\tilde Z_0 )\\
&-2(1-\eta_1) \nabla \eta_2 \nabla (\hat z_0-\hat z_1)-(1-\eta_1)\Delta \eta_2 (\hat z_0-\hat z_1)\|_\star  \leq C|\ln(\la\mu_\la^2)|^{-1}.
\end{align*}
By combining this with \eqref{ee1}, we conclude that
\begin{equation}
\label{2707e1}
\|\tilde L (\tilde z_0)\|_{\star}\leq C |\ln (\lambda\mu_\la^2 )|^{-1}.
\end{equation}

\medskip
Finally, we estimate $\langle \tilde L(\tilde z_0),\tilde z_0\rangle$. Arguing as above, it is easy to see that 
$$
\langle \tilde L(\tilde z_0),\tilde z_0\rangle= \int_{B_{R+1}\backslash B_R} \tilde L(\tilde z_0)\tilde z_0 dx+
\int_{B_{\tau/(\sqrt\lambda\mu_\la )} \backslash B_{\tau/(2\sqrt\lambda\mu_\la )}} \tilde L(\tilde z_0)\tilde z_0 dx+O\left((\lambda\mu_\la^2)^\alpha\right).
$$
Using \eqref{traine1}, \eqref{boundz1}, and \eqref{traine2}, we get
\begin{align}\label{II}
\Bigg| \int_{B_{\tau/(\sqrt\lambda\mu_\la )} \backslash B_{\tau/(2\sqrt\lambda\mu_\la )}} &\tilde  L(\tilde z_0)\tilde z_0 dx \Bigg| \\
&
\leq C \int_{B_{\tau/(\sqrt\lambda\mu_\la )} \backslash B_{\tau/(2\sqrt\lambda\mu_\la )}} \left(|\nabla \eta_2||\nabla \hat z_0||\hat z_0| + |\Delta \eta_2||\hat z_0|^2\right) +O\left((\lambda\mu_\la^2)^\alpha\right) \nonumber\\
&\leq C|\ln (\lambda\mu_\la^2 )|^{-2}.\nonumber
\end{align}
On the other hand, we have 
\begin{multline*}
I:=\int_{B_{R+1}\backslash B_R} \tilde L(\tilde z_0)\tilde z_0 dx=\\
2\int_{B_{R+1}\backslash B_R} \nabla \eta_1 \nabla (\hat z_0-\tilde{Z}_0)\hat z_0 dx+ \int_{B_{R+1}\backslash B_R} \Delta \eta_1 (\hat z_0-\tilde{Z}_0)\hat z_0 dx+O\left((\lambda\mu_\la^2)^\alpha\right).
\end{multline*}
Thus, integrating by parts, we find
$$
I=\int_{B_{R+1}\backslash B_R} \nabla \eta_1 \nabla (\hat z_0- \tilde{Z}_0)\hat{ z}_0 dx -\int_{B_{R+1}\backslash B_R} \nabla \eta_1 (\hat z_0-\tilde{Z}_0)\nabla \hat z_0  dx+ O\left((\lambda\mu_\la^2)^\alpha\right).
$$
Now, we observe that, for $s\in(R,R+1)$, $|\tilde Z_0(s)-\hat z_0(s)|\leq C |\ln (\lambda\mu_\la^2 )|^{-1}$, while $|\hat z_0'(s)|\leq C\left(\frac{1}{R^3}+\frac{1}{R}|\ln (\lambda\mu_\la^2 )|^{-1}\right)$. Thus
$$
\left|\int_{B_{R+1}\backslash B_R } \nabla \eta_1 (\hat z_0-\tilde Z_0)\nabla \tilde z_0 dx \right|\leq \frac{D}{R^3}|\ln (\lambda\mu_\la^2 )|^{-1},
$$
where $D$ is a constant that does not depend on $R$. Note that 
\begin{align*}
\int_{B_{R+1}\backslash B_R} &\nabla \eta_1 \nabla (\hat z_0 - \tilde Z_0)\hat z_0 dx \\
&= 2\pi
\int_R^{R+1}\eta_1'(\hat z_0 - \tilde Z_0)'\left(\tilde Z_0 +O(|\ln (\lambda\mu_\la^2 )|^{-1})\right) s ds\\
&=-  \dfrac{2\pi}{\displaystyle \int_R^{\tau/(\sqrt\lambda\mu_\la )}\dfrac{dt}{ t\tilde{Z}_0^2(t)}}\int_R^{R+1}
\eta_1'(s) \left(\tilde{Z}_0 (s) \int_R^s \frac{dz}{z \tilde{Z}_0^2 (z)} + \frac{1}{s} \right)ds\left[1+O(|\ln (\lambda\mu_\la^2 )|^{-1})\right]\\
&=-\frac{E}{|\ln (\lambda\mu_\la^2 )|}\left[1+o_R(1)+O(|\ln (\lambda\mu_\la^2 )|^{-1} ) \right],
\end{align*}
where $E$ is a strictly positive constant independent of $R$ and $\la$. We thus conclude, choosing $R$ large enough, that $I\sim -E|\ln (\lambda\mu_\la^2 ) |^{-1}$. Combining this with \eqref{II}, we find
$$
\langle \tilde L (\tilde z_0),\tilde z_0\rangle =-\frac{E}{|\ln (\lambda\mu_\la^2 )|}\left[1+O(|\ln (\lambda\mu_\la^2 )|^{-1} ) \right] .
$$
This together with \eqref{estimateA} and \eqref{2707e1}, yields
$$
|A|\leq C|\ln (\lambda\mu_\la^2 )|\|\tilde h\|_{\star}.
$$
Recalling \eqref{modif} and using \eqref{ModifiedEst}, we then deduce that
$$
\|\phi\|_{L^\infty (B_1)} \leq C(\|\tilde h\|_{**} +|\ln (\lambda\mu_\la^2 )|\|\tilde h\|_\star ).
$$
Observe that 
$$
\|\tilde h\|_\star =\sup_{s\in [0,1/(2\sqrt\lambda\mu_\la )]}\frac{\tilde h(s)}{\lambda\mu_\la^2 +(1+s)^{-2-\nu}}\leq \sup_{r\in[0,1/2]} \frac{\lambda\mu_\la^2 |h(r)|}{\lambda\mu_\la^2+\left(1+\frac{r}{\sqrt \lambda\mu_\la}\right)^{-2-\nu}}\leq\|\tilde \chi_1h\|_{\star}.
$$
The previous two inequalities then yield  
$$
\|\phi\|_{L^\infty (B_1)} \leq C(\|h\|_{**} +|\ln (\lambda\mu_\la^2 )|\|\chi_1 h\|_\star).
$$
Recalling the definition of the norm $\| \cdot \|_{*}$, we conclude that 
$$
\|\phi\|_{L^\infty (B_1)}  \leq C\|h\|_*.
$$
It only remains to prove the existence part of the statement. For this purpose, we consider the space 
$$
H=\left\{ \phi \in H^1(B_1)\ | \ \phi \mbox{ is radial} \right\},
$$
endowed with the inner product $\langle \phi,\psi \rangle_{H^1} = \int_{B_1}\nabla \phi \nabla \psi dx+\int_{B_1}\phi \psi dx$. Problem \eqref{lineari0} expressed in weak form is equivalent to finding $\phi\in H$ such that
$$
\langle \phi ,\psi \rangle_{H^1} = \int_{B_1}[\lambda e^U\phi +h]\psi dx\quad \mbox{for all }\psi \in H.
$$
By Fredholm's alternative this is equivalent to the uniqueness of solutions to this problem, which is guaranteed by \eqref{estimateLinear}.
\end{proof}

\section{The proof of Theorem \ref{mainthm}}\label{secMain}
\begin{proof}
Thanks to Proposition \ref{teolineal}, we know that the operator $L$ is invertible. Therefore, we can rewrite \eqref{eqdephi} as
$$\phi=T (\phi)=L^{-1} [R (U)+N (\phi)].$$
Let $\rho$ be a fixed number. We define
$$A_\rho =\left\{\phi \in L^\infty (B_1) :\ \left\|\phi\right\|_{L^\infty (B_1)}\leq \rho \varepsilon_\la^{1+\sigma} \right\},$$
where $\sigma$ is the constant defined in Lemma \ref{estR}. We will show that the map $T : A_\rho \rightarrow A_\rho$ is a contraction. 
Using Lemma \ref{estN}, recalling the definition of $\| \cdot \|_\ast$ given in \eqref{normstar}, and since $|\ln (\lambda\mu_\la^2 )| = O\left( \varepsilon_\la^{-1}\right)$, we see that
\begin{align*}
\left\|\lambda e^{U}\right\|_{\ast}&\leq C \max \left(|\ln (\lambda\mu_\la^2 )|\sup_{r \leq \delta} f_\lambda (r)\dfrac{1}{\lambda\mu_\la^2 \left(1+ \left(\dfrac{r}{\sqrt{\lambda}\mu_\la }\right)^2\right)^2},\right. \\
& \left.  \hspace*{7cm} |\ln (\lambda\mu_\la^2 )|\sup_{\delta \leq r\leq 3/4} f_\lambda (r ) \varepsilon_\la^\beta  , \varepsilon_\la^{-1} \right)\\
&\leq C \varepsilon_\la^{-1}.
\end{align*}
 From this and recalling the definition of $N(\cdot)$ (see \eqref{defN}), we deduce that, for $\phi,\psi \in A_\rho$,
\begin{equation}
\label{mainthme1}
\left\|N (\phi)\right\|_{\ast}\leq \left\|\lambda e^{U}\right\|_{\ast} \left\|\phi\right\|_{L^\infty (B_1)}^2\leq C \varepsilon_\la^{-1}\left\|\phi\right\|_{L^\infty (B_1)}^2
\end{equation}
and
$$\left\|N (\phi)-N (\psi) \right\|_{\ast}\leq  C \varepsilon_\la^{-1}\max \left\{\left\|\phi\right\|_{L^\infty (B_1)},\left\|\psi\right\|_{L^\infty (B_1)} \right\} \left\|\phi-\psi\right\|_{L^\infty (B_1)}. $$
Next, using Lemma \ref{estR}, we obtain
\begin{align}
\label{mainthme2}
\left\|R (U)\right\|_{\ast}& \leq C \max \left(|\ln (\lambda\mu_\la^2 )|\sup_{r \leq \delta}   f_\lambda (r)\dfrac{\dfrac{r^2}{\ep_\la}+r^2|\ln r|+(\la\mu_\la^2)^\alpha}{\la\mu_\la^2\left(1+ \left(\dfrac{r}{\sqrt{\lambda}\mu_\la}\right)^2\right)^2}\right. , \\
&  \left. \hspace*{7cm}|\ln (\lambda\mu_\la^2 )|\sup_{\delta \leq r\leq 3/4} f_\lambda (r ) \varepsilon_\la^\beta, \varepsilon_\la^{1+\sigma} \right)\nonumber\\
&\leq C \varepsilon_\la^{1+\sigma}\nonumber.
\end{align}
%Using Propositions \ref{errorbound1} and \ref{errorin1},
Thus, using \eqref{mainthme1} and \eqref{mainthme2}, we get that, for $\phi,\psi\in A_\rho$,
$$\left\|T (\phi)\right\|_{L^\infty (B_1)} \leq C (\left\|N (\phi)\right\|_{\ast}+\left\|R (U)\right\|_{\ast} )\leq C \varepsilon_\la^{1+\sigma}(\rho^2\varepsilon_\la^\sigma+1)$$
and
$$\left\|T (\phi)-T (\psi)\right\|_{L^\infty (B_1)} \leq C \left\|N (\phi)-N (\psi)\right\|_{\ast}\leq C \rho\varepsilon_\la^\sigma \left\|\phi-\psi\right\|_{L^\infty (B_1 )},$$
where $C$ is a constant independent of $\rho$. It follows that, for any sufficiently small $\la$ (and thus $\ep_\la$), $T$ is a contraction mapping in $A_\rho$, and it therefore has a unique fixed point in $A_\rho$ for $\rho>2C$. This concludes the proof.
\end{proof}

\appendix
\section{An elliptic estimate}
We show a very rough elliptic estimate which is needed in the proof of Lemma \ref{expHo}.
\begin{lem}
\label{appregell}
Let $R>0$ and  $u\in H^1 (B_R (0))$ be a radial solution to
$$
\left\{
\begin{array}{rcll}-\Delta u +u &=& f& \mathrm{in}\ B_R(0)\\ u^\prime (R)&=&g,& \end{array}\right.
$$
for some $f\in L^q (B_R (0))$, with $q>2$. Then, we have
$$\|u\|_{L^\infty (B_R (0))} \leq  C\left[ \left(\dfrac{1}{R}+R\right) R^{1-2/q} \|f\|_{L^q (B_R (0))}+\left( \dfrac{ 1}{R}+ R^2  \right) \|g\|_{L^\infty (\partial B_R (0))} \right] $$
and
$$\|u' \|_{L^\infty (B_R (0))}\leq  C\left[ R^{1-2/q} \|f\|_{L^q (B_R (0))}+(  1+ R  ) \|g\|_{L^\infty (\partial B_R (0))} \right]$$
for some constant $C$ not depending on $R$.
\end{lem}
\begin{proof}
Multiplying the equation by $u$ and integrating by parts, we get
\begin{equation}
\label{regelle1}
\|u\|_{H^1 (B_R )}^2 \leq  \|f\|_{L^2 (B_R )} \|u\|_{H^1 (B_R )}+R |u' (R)| |u(R)| .
\end{equation}
Since $u(R)-u(r)=\int_r^R u' (s) ds$, one deduces that
$$|u(R)|^2 \leq C \left[|u(r)|^2 +\|u'\|_{L^2 (B_R )}^2 \ln \dfrac{R}{r}  \right],$$
where throughout the proof $C$ denotes a constant not depending on $R$. Multiplying the previous inequality by $r$ and integrating, we find
$$R^2 |u(R)|^2 \leq C[\| u\|_{L^2 (B_R )}^2+\|u'\|_{L^2 (B_R) }^2 R^2  ].$$
This implies that
\begin{equation}
\label{regelle2}
|u(R)|\leq  C\left( \dfrac{1}{R}+1 \right) \|u\|_{H_1 (B_R )}.
\end{equation}
From \eqref{regelle1}, \eqref{regelle2}, and $u' (R)=g$, we obtain that
$$
\|u\|_{H^1 (B_R )}^2 \leq  \|f\|_{L^2 (B_R )} \|u\|_{H^1 (B_R)}+C(  1+ R  ) \|g\|_{L^\infty (\partial B_R )} \|u \|_{H^1 (B_R)} .$$
 Thanks to H\" older inequality, we find that
\begin{equation}
\label{lastestregH1}\|u\|_{H^1 (B_R)} \leq C[ R^{1-2/q} \|f\|_{L^q (B_R)} +(1+R) \|g\|_{L^\infty (\partial B_R)} ]  .\end{equation}
Next, observe that for any $s\in (0,R)$ we can rewrite the equation as
$$u' (s) s= \int_0^s (u-f) r dr.$$
From H\" older inequality, we obtain that
$$|u'(s)|\leq C \|u-f \|_{L^2 (B_R )  } \leq C (\|u\|_{L^2 (B_R )}+ R^{1-2/q} \|f\|_{L^q (B_R)} ).$$
From \eqref{lastestregH1}, we deduce that
$$\|u'\|_{L^\infty (B_R )}\leq C (R^{1-2/q} \|f\|_{L^q (B_R )} +(1+R )  \|g\|_{L^\infty (\partial B_R )}  ).$$
By noting that
$$u(R)-u(\tilde s)= \int_{\tilde s}^R u' (r) dr,$$
we get from \eqref{regelle2} that
\begin{align*}
\|u\|_{L^\infty (B_R)}& \leq C\left[ \left( \dfrac{1}{R}+1 \right) \|u\|_{H^1 (B_R)}+ R \|u'\|_{L^\infty (B_R)}\right] \\
&\leq C \left[ \left(\dfrac{1}{R}+1+R\right)  R^{1-2/q} \|f\|_{L^q (B_R)}  + \left(\dfrac{1}{R}+R^2\right) \|g\|_{L^\infty (\partial B_R)}\right].
\end{align*}
This concludes the proof.
\end{proof}

\bibliography{ReferencesKellerSegel}
\end{document}